\documentstyle{article}
\title{
Singular hermitian metrics on vector bundles  }
\author{Mark Andrea A.  de Cataldo\thanks{Partially 
supported by N.S.F. Grant DMS-9701779}
 }
\date{}

\newtheorem{tm}{Theorem}[subsection]
\newtheorem{lm}[tm]{Lemma}
\newtheorem{pr}[tm]{Proposition}
\newtheorem{rmk}[tm]{Remark}

\newtheorem{ex}[tm]{Example}
\newtheorem{defi}[tm]{Definition}

\newtheorem{??}[tm]{Question}

\font\tenmsb=msbm10
\font\sevenmsb=msbm7
\font\fivemsb=msbm5

\newfam\msbfam
\textfont\msbfam=\tenmsb
\scriptfont\msbfam=\sevenmsb
\scriptscriptfont\msbfam=\fivemsb
\def\Bbb#1{{\fam\msbfam #1}}

\font\teneufm=eufm10
\font\seveneufm=eufm7
\font\fiveeufm=eufm5
\newfam\eufmfam
\textfont\eufmfam=\teneufm
\scriptfont\eufmfam=\seveneufm
\scriptscriptfont\eufmfam=\fiveeufm
\def\frak#1{{\fam\eufmfam\relax#1}}

\newcommand\ci{\cite}
\newcommand\s{\sigma}

\newcommand\comp{{\Bbb C}}

\newcommand\zed{{\Bbb Z}}
\newcommand\nat{{\Bbb N}}
\newcommand\pn[1]{{\Bbb P}^{#1}}

\newcommand\blacksquare{{\hspace*{\fill} $\Box$}} 
\newcommand\odix[1]{ {\cal O}_{#1} }
\newcommand\odixl[2]{ {\cal O}_{#1}({#2}) }

\newcommand\T[2]{\Theta_{#1}({#2})}
\newcommand\TT[2]{\tilde{\Theta}_{#1}({#2})}

\newcommand\D[1]{d'd'' {#1}}
\newcommand\e{\epsilon}
\newcommand\id{{\cal I}}

\newcommand\E{{\cal E}}
\newcommand{\SS}{\Sigma}

\newcommand{\LL}{\frak{L}}

\newcommand\surj{-\!\!\!-\!\!\!-\!\!\!-\!\!\!\!\!\gg}
\newcommand\QQ{\frak Q}
\newcommand\q{\frak q}

\begin{document}
\maketitle
\begin{center}
{\em Dedicated to the memory of Michael Schneider}
\end{center}
\bigskip

\begin{abstract}
We introduce a
 notion of singular hermitian metrics (s.h.m.) for holomorphic vector bundles
and define positivity in view of $L^2$-estimates.
 Associated with a suitably positive  s.h.m.
there is a (coherent) sheaf $0$-th kernel
of a certain $d''$-complex. We prove  a vanishing theorem for the 
cohomology 
of this sheaf.
All this generalizes  to the case of higher rank known results of Nadel
for the case of line bundles. 
We introduce a new semi-positivity notion,
$t$-nefness, for vector bundles, establish some of its 
basic properties and prove that on curves  it coincides  with ordinary 
nefness.
We particularize the results  on s.h.m. to the case of vector bundles
of the form $E=F \otimes L$, where $F$ is a $t$-nef vector bundle and $L$ 
is a 
positive (in the sense of currents) line bundle. As applications we 
generalize
to the higher rank case 
1) Kawamata-Viehweg Vanishing Theorem, 2) the effective results
concerning the global generation of jets for the adjoint to
 powers of ample line bundles, and
 3) Matsusaka Big Theorem made effective.

\end{abstract}

\tableofcontents

\section{Introduction}
In this study I introduce a notion of singular hermitian metrics
({\em s.h.m.}) on holomorphic vector bundles over complex manifolds.
The original motivation was to explore the possibility of employing,
in the setting of vector bundles,
the new transcendental techniques developed by Demailly  and Siu
in order to  study
global generation problems for (adjoint) line bundles.
The notes
\ci{dem94} are an excellent introduction to these techniques and to the 
results
in the literature.  
One can consult the lucid notes \ci{eln}  
for an algebraic counterpart to these techniques.

\bigskip
Let me discuss the case of line bundles.
Let $X$ be  a non-singular projective manifold of dimension $n$, $ L$
and $E$ be 
an ample  and a nef line bundle  on $X$, respectively,
 $a$   be  a non-negative  integer and $m$  be a positive one.

\noindent
{\em 
Under which conditions on $a$ and $m$ will the line bundle
$$
{\frak P}: =K_X^{\otimes a }\otimes { L}^{\otimes m} \otimes {E}
$$
be generated by its global sections (free)?} 

\noindent
More generally, we can ask for conditions
on $a$ and $m$ under which
 the simultaneous generation of the higher jets of $\frak P$ 
at a prescribed number of  points on $X$
is ensured.

It is clear that $m\gg 0$ answers the question. However,
how big $m$ should be could depend, {\em a priori}, on $X$. 
For example, 
Matsusaka Big Theorem asserts that $L^{\otimes m}$ is   very ample
for every $ m \geq M:=M(n, L^n, K_X\cdot L^{n-1})$.  An effective value
for $M$  has been recently determined in \ci{siumat} and \ci{dem96};  see 
also
\ci{fdb} for the case of surfaces.

The presence of the canonical line bundle,  i.e. $a>0$, changes dramatically
the shape of the lower bound on $m$.  
Fujita's Conjecture speculates that $K_X\otimes {L}^{\otimes m}$ should
be free
 as soon as $ m\geq  n+1$.
This  conjecture
is true for $n\leq 4$ by the work of Reider, Ein-Lazarsfeld and Kawamata. 
In the papers  \ci{an-siu} and \ci{tsu} it is proved that 
$m\geq \frac{1}{2}(n^2 + n +2)$ gives freeness. 

Effective results depending only on $n$ are proved for $a\geq 1$ by 
several authors.
The seminal paper is \ci{dem93b}
where it is proved, by (differential-geometric-){\em analytic methods}, 
that
$K_X^{\otimes 2}\otimes L^{\otimes m}$ is very ample for all $ m \geq 12n^n$.
Then followed the  paper \ci{ko}, where a similar result
is proved using {\em algebraic-geometric  methods}. Since then, several  
papers
have appeared on the subject. 
The  reader may  consult the following references to compare the various
results and techniques:
\ci{dem94} (an account of the analytic approach with  a rather complete
bibliography), \ci{eln} (an account of the algebraic approach
and of how many of the analytic instruments may be re-tooled
and made algebraic), \ci{an-siu} and \ci{tsu}
(freeness;  written in the  analytic language, but apt  to be completely
translated into the algebraic language after observations by Koll\'ar 
\ci{kos},
\S5;
see also \ci{siu-tak}), \ci{dem96},  \ci{siu-va} and \ci{siu94b} 
(very ampleness;
analytic), \ci{siumat} and \ci{dem96} (an effective version of Matsusaka
Big Theorem; analytic).

\bigskip
An extra nef factor $E$ plays a minor role and 
all of the results quoted above hold
 in its presence. 
This simple
fact was the starting point
of my investigation. 

\begin{??}
\label{qk}
 Can we obtain effective results on $a$ and $m$ for 
the global generation of the vector bundles $\frak P$ by assuming that
$E$ is a  suitably semi-positive vector bundle of rank $r$?
More generally, can we obtain similar results about 
the simultaneous generation of the higher jets of $\frak P$ 
at a prescribed number of  points on $X$\,?
\end{??}

\medskip 
I expected that the statements in the aforementioned literature concerning
 the line bundles $\frak P$ with the  
nef line bundle  $E$
should carry over, 
{\em unchanged}, to the case in which $E$ is a nef vector bundle.

\noindent
On a projective  manifold a
 nef line bundle can be endowed with  hermitian metrics whose curvature forms
 can be made to have  arbitrarily
small negative parts
(cf. Definition
\ref{defneflb}). In the analytic context this  fact can be used to make 
the 
presence  of  a nef line bundle $E$ harmless. The same is true in 
the algebraic context
because of the 
the  numerical properties of 
 nefness. 

\noindent
 A natural  algebraic approach to the case of higher rank
is to consider an analogous question for the tautological line bundle $\xi$
of the projectivized bundle $\pi: {\Bbb P}({E}) \to X$.
The  results I obtain
with the algebraic approach  are for
${\frak P}\otimes \det E$; 
 compare  Remark  \ref{algnef}  with the sample 
effective global generation result presented below;

 see \ci{deeff}.
\noindent
On the analytic side, the problem is that the nefness
 of a vector bundle $E$ 
  does not seem to be  linked to a curvature
condition on $E$ itself.

\smallskip
As far as Question \ref{qk}
is concerned, nefness does not seem to give  enough room to work  
analytically
with higher 
rank vector bundles.

\noindent 
Instead I introduce, for every  vector bundle $E$ and
every positive integer $t$, the notion of 
{\em $t$-nefness}  which is a  new  semi-positivity concept for
 vector bundles.
In some sense it is in between the algebraic
notion of nefness and the differential-geometric notion
of  $t$-semipositivity. It is a natural higher rank
 curvature analogue of  the aforementioned
characterization of nef line bundles.
The property of $t$-nefness is checked by considering tensors 
 in $T_X\otimes {E}$ of rank at most $t$; such tensors have ranks
never bigger
than $N:=\min (\dim X, {\rm rank} \,{E})$. Incidentally,
$(t+1)$-nefness implies $t$-nefness for all positive integers $t$, 
$1$-nefness 
implies nefness and I do not know whether nefness implies $1$-nefness.

\smallskip
Though, as I show in Theorem \ref{umemura},  on curves $1$-nefness is
equivalent to nefness, the  notion of $t$-nefness
 is rather difficult to check in an algebraic 
context. 
However, see Example \ref{listnef} for a list of   nef bundles 
which I know to be  $N$-nef or from which it is easy
to obtain $N$-nef bundles (e.g. nef bundles on curves, nef line bundles,
flat bundles,
nef bundles
on toric or abelian varieties, the 
tangent bundles of low-dimensional
K\"ahler manifolds with nef tangent bundles,  pull-backs, etc.).

\smallskip
Assuming that $E$ is
 $N$-nef, I   
prove for  the vector bundles $\frak P$
the same statements as the ones in the literature for  the line bundle 
case; 
see Theorem \ref{effres}.
Moreover, if $\,{E}$ is $1$-nef, then 
the same results hold replacing $E$ by  $E \otimes \det  E$.
   The scheme of the proofs 
is the same
as in the rank one case (see Proposition \ref{cucu}, 
\S\ref{vecbl}, and of course \ci{dem94}, \S5 and \S8). 
However,  at each and every step we need higher rank analogues of the
analytic package developed for the line bundle case by
Demailly and  Nadel: 
regularization, $L^2$-estimates, coherence of relevant sheaves and vanishing
theorems.
For the purpose of proving  these effective results for the vector bundles
$\frak P$,    
one would have to make precise the notion
of singular hermitian metrics with positivity  
and  prove their relevant properties 
in a  special case: the one of a hermitian vector bundle twisted by a 
line bundle
endowed with a singular metric.  Then one would have to   
prove the relevant vanishing theorems. All this  can be done by building on
\ci{dem82}, \S5 and \S9.

\noindent
 However,
 I felt 
that  it should be worthwhile to develop 
a general theory of {\em singular hermitian metrics} on vector bundles 
with special 
regards towards positivity.
 
\medskip
Inspired by the case of line bundles, in this paper I develop such a theory
 and obtain as
an application the effective results  mentioned above.
To get a flavor of the results let me state
(\ref{effres}.$1'$), which constitutes an answer to Question \ref{qk}
(see Remark \ref{geoint} for a geometric interpretation of 
these kind of results):

\medskip
\noindent
{\bf Effective global generation.} 
{\em 
Let $E$ be $N$-nef. Then $K_X \otimes L^{\otimes m}
\otimes {E}$ is globally generated by 
its global sections
for all $m\geq \frac{1}{2}(n^2 + n +2)$. Moreover, if $E$ is $1$-nef,
then the previous statement is true if we replace $E$ 
by ${E} \otimes \det  E$.
}

\bigskip
The paper is organized as follows.

\noindent
\S{\bf  1}  fixes the notation. \S{\bf 2} is devoted to s.h.m. which are 
defined in
\S2.1. The case of line bundles is discussed in \S\ref{exlbshm}. In 
\S2.3 we introduce
the sheaf  $\E (h)$ which generalizes Nadel multiplier ideal sheaf.
In \S2.4 we define positivity for s.h.m. and study some of its properties.
 \S{\bf 3} revolves about
the notion of
$t$-nefness. The definition  and the basic properties are to be found in
\S3.1  and \S3.2, respectively. \S3.3 is devoted to the proof
of Theorem  \ref{umemura}  which ensures that on curves the 
algebraic-geometric
 notion
 of nefness
can be characterized
 by the differential-geometric notion of $1$-nefness.
\S3.4 consists of a  footnote to \ci{d-p-s}, Theorem 1.12: ampleness
for  a vector bundle $E$ can be characterized by a curvature condition
on a system of metrics on {\em all} symmetric powers $S^pE$ of $E$, 
though positivity
may occur only for $p\gg 0$.
\S{\bf 4} is devoted to vanishing theorems. The basic one is Theorem 
\ref{vanish},
a generalization of Nadel Vanishing Theorem; 
Proposition \ref{coherent} asserts
that $\E(h)$  is coherent  in the presence of
suitable positivity.
 \S4.2 links   $t$-nefness and positivity via  vanishing; see
Theorem \ref{myvan}. Theorem \ref{kv} is a  generalization of
Kawamata-Viehweg Vanishing Theorem.
\S{\bf 5} contains the effective results concerning the vector bundles
$\frak P$. \S5.1 contains, 
for the reader's convenience,
 a summary of the results of Anghern-Siu and Siu concerning  
special s.h.m. on line bundles which,
transplanted to $N$-nef
vector bundles, 
will provide the global generation of jets. We also offer  the simple Lemma
\ref{freetojet}, which constructs metrics with similar properties 
starting from free 
line bundles. \S\ref{vecbl} contains our effective results
concerning the vector bundles $\frak P$; see Theorem
\ref{effres}.

\medskip
\noindent
{\bf Acknowledgments}. I heartily thank
J.-P. Demailly for
reading a preliminary and rough version of
this paper and for suggesting some improvements. 
I am indebted to J. Koll\'ar
for posing a   question similar to Question \ref{qk}.
I thank  L. Ein and R. Lazarsfeld for convincing me to
think about an algebraic proof of the results of  Theorem \ref{effres};
this has lead me to the statements of Remark \ref{algnef}; see \ci{deeff}.
It is a pleasure to thank the participants of 
the lively algebraic geometry seminar at Washington University in St. Louis
for their encouragment and useful criticisms:  
V. Masek, T. Nguyen, P. Rao
and D. Wright. 
I would like to thank N.M. Kumar for 
many pleasant and useful conversations.

\section{Notation and preliminaries}
Our basic reference for the language 
of complex differential geometry is \ci{g-h}. Sufficient and
 more self-contained references  are \ci{dem82}, \S 2 and \ci{dem94},  \S 3.

All manifolds are second countable, connected and  complex; 
the dimension is
the complex one.
 All vector bundles are holomorphic.
The term {\em hermitian metric}
always refers to a  hermitian metric  of class ${\cal C}^2$. 
  A {\em hermitian bundle} $(E,h)$ is the assignment of a vector bundle 
$E$ together
with a hermitian metric $h$ on it.

Duality for vector bundles is denoted by the symbol $`` \, ^ * \, "$ 
and ${\rm End} (E)$
is the vector bundle of endomorphisms of $E$. We often do not distinguish
between vector bundles and associated sheaves of holomorphic sections; 
at times, 
we employ simultaneously the additive and  multiplicative notation for 
line bundles.

\smallskip
- $d=d' + d''$ denotes the natural decomposition
of the exterior derivative $d$ into its  $(1,0)$ and $(0,1)$ parts; $d''$
denotes also the usual operator associated with a vector bundle $E$.

\smallskip
If $(E,h)$ is a rank $r$ hermitian 
vector bundle on a manifold $X$ of dimension $n$, then we denote by:

- $D_h(E)$ the associated {\em hermitian connection} which is  also  called
the {\em Chern connection};

- $\T{h}{E}=iD^2_h(E)$ 
 the associated {\em curvature tensor};

\noindent
in particular, if $L$ is a line bundle with a metric $h$, represented 
locally
on some open set $U$ by $e^{-2\varphi}$, then we have
$ \T{h}{L}_{|U}= 2i \D{\varphi}$;

- $\tilde{\Theta}_h (E)$ the associated hermitian form on $T_X \otimes E$.

\smallskip
 If $\theta$ is a hermitian form on a complex vector space $V$, we denote
$\theta (v,v)$ by $\theta (v)$; if in addition, $\theta$ is positive 
definite,
then
we denote $\theta(v)$ by $|v|^2_{\theta}$.

- ${\rm Herm}_h (V)$ is  the set of endomorphisms $\alpha$  
of a hermitian vector space $(V,h)$ such that $h(\alpha(v),w)=h(v,\alpha 
(w))$, 
$\forall
v,\, w \in V$.

\noindent
Given $\Theta$, a real $(1,1)$-form with values in ${\rm  Herm}_h (E)$, 
we denote
the associated hermitian form on $T_X \otimes E$ by $\Theta_h$, or by
$\Theta$, if no confusion is likely to arise.  
The hermitian form
$\TT{h}{E}$  will be denoted from now on by
$\T{h}{E}$.  If $\omega$ is  a real $(1,1)$-form,  e.g. the one associated
with a hermitian metric on $X$, then  $\omega \otimes {\rm Id}_E$ has values
in ${\rm Herm}_h(E)$ and we denote the associated hermitian form by
 $\omega \otimes {\rm Id}_{E_h}$ so that    
$\omega \otimes {\rm Id}_{E_h}(t\otimes e)=$
$\omega (t,it) |e|^2_h$, $\forall x \in X$, $\forall
t \in T_{X,x}$ and $\forall e \in E_x$.

\medskip
\noindent
{\bf The rank of a tensor.}
Let  $V$ and $W$ be complex vector spaces of finite 
dimensions  $r$ and $s$, 
respectively,  
$v=\{v_{i}\}_{i=1}^r$ and
$w=\{w_{\alpha}\}_{\alpha=1}^s$
 be bases for  $V$ and $W$, respectively;
 tensor products are taken over $\comp$.

\noindent
Every tensor $\tau \in V \otimes W$ defines two linear maps 
$\alpha_{\tau}: W^* \to V$ and 
$\beta_{\tau}: V^* \to W$; moreover,  we can write 
$\tau=\sum_{i\alpha} \tau_{i\alpha}v_i\otimes w_{\alpha}$ and
associate with $\tau$ the 
$r\times s$ matrix $||\tau_{i\alpha}||$.  
The integer 
  $\rho (\tau):= {\rm rank} (\alpha_{\tau})=$ ${\rm rank} 
(\beta_{\tau})=$ 
${\rm rank} ||\tau_{i\alpha}||$ is called the {\em
rank} of the tensor $\tau$. 

\noindent
Tensors of rank zero or one are called 
{\em decomposable}; they have the form
$\tau=v\otimes w$, for some $v \in V$ and $w\in W$. For any non-zero
tensor $\tau\in V \otimes W$ 
we have
that  $1 \leq \rho (\tau) \leq  \min (r, s)$. In particular, if either
$r=1$, $s =1$, or both, then every tensor $\tau\in V \otimes W$ is 
decomposable.

\medskip
\noindent
{\bf Inequalities associated with the rank.}
Given two hermitian forms $\theta_1$ and $\theta_2$ on  
$V\otimes W$,  we can 
compare them on tensors of various rank. Let $t$ be any positive integer.
 We write  $\theta_1 \geq_t  \theta_2$ 
if
the hermitian form $\theta_1 - \theta_2$ is semi-positive 
 definite on all tensors in $V \otimes W$
of rank $\rho \leq t$. 
If $\theta \geq_t 0$, 
then $\theta \geq_{t'} 0$  for every $ t' \leq t$.
If $\theta_1 \geq_{\min (r,s) } \theta_2$, then 
$\theta_1 \geq _t  \theta_2$ for every $t$.
The symbol $>_t$ can be defined analogously and it enjoys
similar properties.

\medskip
These considerations  and this language are easily transferred
to vector bundles.

\section{Singular hermitian  metrics on  vector bundles}
\label{shm}
In this section we  define singular hermitian metrics on vector bundles,
discuss the  case of line bundles, introduce the sheaf $\E (h)$
and define
 positivity.

\subsection{The definition of singular hermitian metrics}
\label{dshm}
Let $X$ be a  manifold of dimension $n$, $E$  be a rank $r$ 
vector bundle over $X$ and $\bar{E}$ the conjugate of $E$.
Let $h$ be a section of the smooth 
vector bundle $E^*\otimes \bar{E}^*$ with measurable coefficients,
such that $h$ is an almost everywhere (a.e.)
positive definite hermitian form  on $E$; we call such an 
$h$ a
{\em  measurable metric} on $E$.
A measurable metric  $h$ on $E$ induces naturally measurable metrics
on $E^*$, on any tensor representation of $E$,  e.g.
$T^{\alpha} E$, $S^{\beta} E$, $\wedge^{\gamma} E$ etc., on any quotient 
bundle 
of $E$,  etc.

In practice these metrics  $h$
occur as {\em degenerate metrics} of some sorts, 
e.g. $h$ is a hermitian metric outside a proper analytic subset $\SS$ of $X$,
so that the curvature tensor is well-defined outside $\SS$.

\smallskip
 We are interested in those $h$ 
for which 
the  curvature tensor has a global meaning.  We propose the following
simple-minded definition.

\begin{defi}
\label{defshm}
{\rm(s.h.m.)}
Let $X$, $E$  and   $h$ be as above and $\SS \subseteq X$ be a closed set
of measure zero.
 Assume that there exists a sequence of hermitian metrics $h_s$
such that:
$$
\lim_{s \to \infty} h_s = h \qquad  in \,\,  the \, \, 
{\cal C}^{2}-topology
\,\,   on\,\,   X \setminus \SS. 
$$ 
 We call the collection of data $(X,E,\SS, h,h_s)$ a singular hermitian 
metric
(s.h.m.) on $E$. 
We call $\T{h}{E_{|X\setminus \SS}}$ the curvature tensor of 
$(X,E,\SS,h,h_s)$ and we
 denote it by $\T{h}{E}$.  $\T{h}{E}$ has continuous coefficients
and values in
${\rm Herm}_h(E)$ away from $\SS$; we denote the a.e.-defined associated 
hermitian 
form  on $T_X \otimes E$  by the same symbol $\T{h}{E}$.
\end{defi}

If no confusion is likely to arise,
we indicate a s.h.m  by $(E,h)$ or simply by $h$.

\medskip
The guiding principle which subtends this definition can be formulated as 
follows.

\smallskip
\noindent
{\em Assume that we would like to
prove a property 
$P$ for $h$ which is true
for all metrics $h'$ of class ${\cal C}^2$ 
in the presence of a certain curvature condition
$C$ on $h'$; 
if  $h$ has the required property
$C$ and we can find hermitian metrics  $h_s$ which regularize
$h$ ``maintaining" $C$, then  $P$ holds for all $h_s$ and we can try to 
prove, 
using limiting arguments,
that $P$ holds for $h$.} 

\smallskip
\noindent
This principle has been successfully exploited in
\ci{dem82}, 
\S5; see \S\ref{exlbshm}
for a brief discussion. 
 We will take this principle as the definition of positivity; see
Definition \ref{ipo}
 and 
Proposition \ref{l2}, where $P$ is the solution to the $d''$-problem
with $L^2$-estimates and $C$ is ``positivity."

\medskip
Because of the convergence in the ${\cal C}^2$-topology,
the  notion of s.h.m. is well behaved under the operations of
taking quotients,  dualizing,  forming direct
sums,   taking tensor products, 
  forming tensor representations,  etc.

\subsection{Discussion of the line bundle case: curvature current, 
positivity,
Nadel Ideal,  Nadel  Vanishing Theorem, and the production of sections}
\label{exlbshm}
We now remark that the singular metrics on line bundles to be found in 
the literature
are s.h.m. We also discuss some of the relevant features
of these metrics in the presence of positivity.
Basic  references for what follows are
\ci{dem94}, \S5, \ci{dem82}, \S9 and \S5. A technical remark:
for the mere purpose of being consistent with Definition \ref{defshm}, in 
what follows
we assume that   plurisubharmonic (psh) functions
are ${\cal C}^2$ outside  a closed set of measure zero. In all the 
applications
one uses {\em algebraic singular metrics} as in \ci{dem96}, so that this
condition is automatically satisfied. However, all the theory described below
and its applications
work without this  restriction; see also \ci{dem92}, \S3.

\smallskip
Note that in what follows we can replace the hermitian line bundle
$(L,h_0)$ by a hermitian vector bundle
$(E,h_0)$ by operating minor changes.

\medskip
A {\em singular metric}  on a line bundle $L$
over a manifold $X$ is, by definition,  a metric of the form
$h=h_0 e^{-2\varphi}$, where $h_0$ is a hermitian metric 
on $L$ and $\varphi$ is a locally
integrable function  on $X$. 
We shall always assume that 
$X$ is K\"ahler and that
$\varphi$ is {\em almost psh},  i.e.  it can be written, locally on $X$,
 as the sum $\varphi=\alpha + \psi$, where 
$\alpha$  is a local function  of class ${\cal C}^2$ and $\psi$ is a local
  psh function.
By taking $d'd''$ in the sense of distributions,
 we can define the
associated {\em curvature  
$(1,1)$-current}:
$$
T:= \T{h_0}{L} + 2i d'd'' \varphi_{ac} + 2i d'd'' \varphi_{sing},
$$ 
where  $2id'd'' \varphi_{\rm ac}$ and $2i d'd'' \varphi_{sing}$
are the absolutely continuous and singular part of $2i d'd'' \varphi$,
respectively;  
$2id'd'' \varphi_{\rm ac}$ has locally integrable
coefficients  and
 $2id'd''\varphi_{\rm sing}$  is  supported on  some closed
set $\SS$ of measure zero.
A 
regularizing-approximating result of  Demailly's
exhibits these singular metrics on line bundles
as s.h.m. by constructing the necessary
regularizing hermitian metrics $\{h_s \}_{s=1}^{\infty}$. 
We have $\T{h}{L}=\T{h_0}{L} + id'd''\varphi_{ac}$.
Similar considerations hold for metrics dual to metrics as above.

\begin{ex}
\label{ef}
{\rm 
(Cf.
\ci{dem94}, Example 3.11 and \ci{dem96}, page 246)
Let  $D=\sum m_iD_i$ be   a divisor with coefficients 
 $m_i \in\zed$. The associated line 
bundle carries a singular metric  with curvature current
$T=2\pi \sum m_i[D_i] $ where  the  $[D_i]$ are  the currents of 
integration over  the subvarieties
$D_i$. These currents are positive 
if and only if all 
$m_i\geq 0$. More generally, given a finite number
of  non-trivial holomorphic sections
of a multiple of a line bundle $L$, we can construct a s.h.m. on $L$.
This metric will be singular only  at  the common zeroes of the sections 
in question.
}
\end{ex}

\smallskip
 The {\em  Nadel ideal} $\id (h)$ (see 
\S \ref{mi}) {\em is coherent}. 
 This is an essential feature in view of the use of
 this ideal
in conjunction with Riemann-Roch Theorem.
 
\medskip

Let $\omega$ be a K\"ahler metric on a weakly pseudoconvex
manifold $X$. Assume that 
$\T{h}{L} \geq \e\omega$ as a $(1,1)$-current, 
for some positive and continuous function $\e$ on $X$.
Then we have {\em Nadel
Vanishing Theorem}: 
$
H^q(X,K_X\otimes L
\otimes \id (h))=0, \quad \forall q >0;$ see \ci{na}.
This can be seen as a consequence of the solution to the
$d''$-problem for $(L,h)$ with $L^2$-estimates; see \ci{dem94}, \S5.

\medskip
As an easy consequence of Nadel Vanishing Theorem we have the following
result which lays the basis for the effective results for the global 
generation
of adjoint line bundles  etc. See \ci{dem94}, Corollary 5.12.

\begin{pr}
\label{cucu}
Let $(X,\omega)$ be  as above and
$\cal L$ be a line bundle over $X$ equipped with a s.h.m. $h$ such
that 
$\T{h}{L} \geq \e \omega$ for some continuous and positive function $\e$ 
on $X$.
Assume that $p$ is a positive integer and  that
$s_1, \ldots, s_p$ are non-negative ones. Let
 $x_1, \ldots, x_p$ be distinct isolated points of the complex space
$V(\id (h))$ such  that $\id (h) \subseteq {\frak m}_{x_i}^{s_i +1}$.
Then there is a surjective map
$$
H^0 (X, K_X + {\cal L} ) \surj
  \bigoplus_{i=1}^p {\cal O}(K_X +
{\cal L} )\otimes
{\cal O}_{X,x_i}/ {\frak m}_{x_i}^{s_i +1}.
$$ 
\end{pr}
Once the analytic package (definition of s.h.m., 
regularization-approximation,
solution of $d''$ with $L^2$-estimates, coherence
of Nadel ideal and Nadel Vanishing Theorem) has been developed, 
in order to solve the global generation problem one needs s.h.m. as in 
Proposition \ref{cucu}. This requires hard work
and it  has been done by Anghern-Siu, Demailly, Siu
and Tsuji. The coherence and the
 vanishing theorem are utilized together with a clever  use of Noetherian 
Induction.

\medskip
We are about to provide a similar analytic package for the case of vector 
bundles.

\subsection{The subsheaf $\E (h)$ associated with a measurable metric
$(E,h)$}
\label{mi}
If $h$ is a measurable metric on $E$ and   $e$ is a measurable section
of $E$, then   the function $|e|_h$ is  measurable.

\begin{defi}
\label{iande}
Let $h$ be a measurable metric on $E$.

\noindent
 Let $\id (h)$  be the analytic sheaf of germs of 
holomorphic functions
on $X$ defined as follows:
$$
 \id(h)_x:= \{ f_x \in \odix{X,x} \!: \, \, |f_x e_x|^2_{h}\, \mbox{is 
integrable in some
 neighborhood
of \,} \, x, \, \forall \, e_x \in E_x\}.
$$
Analogously, we define an analytic  sheaf $\E (h)$ by setting:
$$
 \E(h)_x := \{ e_x \in  E_x \, : \, \, |e_x|^2_{h}\,
\mbox{ is integrable in some
 neighborhood
of } \, x\,\}. 
$$
\end{defi}

\begin{rmk}
\label{sub}
{\rm 
It is easy to show, using the triangle inequality,  that
$\id (h)\otimes E \subseteq \E(h)$.
}
\end{rmk}
We call
$\id(h)$ the {\em
multiplier ideal} of $(E,h)$. 
Note that
if $E$ is a line bundle together with a measurable metric $h$, then
 $\E(h)=\id (h) \otimes E$.

\medskip
There are other subsheaves of $E$, associated with a measurable metric $h$.

\smallskip
Given any measurable metric $h$  on a vector bundle $E$,
 the tautological line bundle
$\xi := \odixl{{\Bbb P}(E)}{1}$ inherits a 
natural measurable metric ${\frak h}$, the quotient
metric of the surjection $\pi^* E \to \xi$; here $\pi: {\Bbb P}(E) \to X$ 
is the structural 
morphism of the projectivized bundle and we are using Grothendieck's
notation. 
We thus  get two sheaves $\id ({\frak h})$ and $\xi \otimes \id (\frak h)$.
If we apply $\pi_*$, then we  get two other subsheaves of $E$.

\medskip
In summary, associated with $(E,h)$ there 
are four subsheaves of $E$:

\smallskip
\centerline{
$\id (h)\otimes E \subseteq \E(h)$, $ \quad \pi_* \id ({\frak h}) \otimes E
\quad $ and
$\quad \pi_* \,\,  \xi \otimes   \id ({\frak h})$.}

\begin{rmk}
{\rm 
The inclusion above  may  be strict. In fact, consider
the vector bundle $\Delta \times \comp^2$, where $(\Delta,z)$ is
 the unit disk
in $\comp^1$; define a s.h.m. by setting $h={\rm diag} (e^{-2\log{|z|}},
e^{-4\log{|z|}})$. Then one checks that $\id (h) =z^2 \cdot 
\odix{\Delta}$ 
and that $\E (h)=
z\cdot \odix{\Delta} \oplus z^2\cdot \odix{\Delta}$.  
The same example shows that $\E (h)$ is not in general equal to neither
$\pi_* \,\, \xi \otimes   \id ({\frak h})$, nor $\pi_* \id ({\frak h}) 
\otimes E$.
In fact, 
a direct computation shows that:  $\id ({\frak h})=\pi^*(z)$.
We have  
     $ \id(h) \otimes E$ $\subset$
$\E(h)$ $\subset$  
 $\pi_* \,  \id  ({\frak h}) \otimes \xi=$
$\pi_* \,  \id ({\frak h}) \otimes E$.
}
\end{rmk}

What is, among  the  four sheaves above,  the ``right" object to look at? 
To  answer  this question we consider: 

\medskip
\noindent
{\bf 
the complex $({\frak L}^{\bullet}, d'')$}. 
Let $h$  be  a measurable metric  on a vector bundle $E$ and $\omega $ 
be a hermitian metric
 on $X$. By following the standard conventions
in \ci{we},  we obtain  a metric with measurable coefficients
for the fibers of ${T_X^{p,q}}^* \otimes E$; we denote this metric again 
by $h$.
We define a complex
$({\frak L}^{\bullet}, d'')$ of sheaves on $X$ as follows. 
This complex is independent
of the choice of $\omega$.

\noindent
Let ${\frak L}^q$ be the sheaf of germs of $(n,q)$-forms $u$ with values 
in $E$ 
and square-integrable  coefficients
 such that  $|u|_h^2$ is locally integrable, 
$d''u$ is defined  in the sense of distributions with square-integrable
coefficients and $|d''u|_h^2$ is locally integrable.

\medskip
The kernel of $d''$ in degree zero is $K_X \otimes \E(h)$ (cf. \ci{g-h}, 
page 380). 
A solution to the $d''$-problem with $L^2$-estimates for 
$(E,h)$ would imply the vanishing of the higher cohomology of
$K_X \otimes \E(h)$. See Theorem \ref{vanish}.

\medskip
If we are aiming at vanishing theorems as in the line bundle case, the 
sheaf $\E(h)$  seems to be  the right object to look at.

\subsection{Positivity}
\label{p}
As is well-known, 
the curvature tensor  $\T{h}{L}$ of a hermitian line bundle $(L,h)$ is 
decomposable
and can be identified with
a
real $(1,1)$-form on $X$. 
This latter is a positive $(1,1)$-form  if and only if the hermitian form
$\T{h}{L}$  is positive on $T_X \otimes L$.

It is therefore natural to define positivity
for singular metrics on line bundles using the notion of {\em positive 
currents} 
according to 
Lelong; see \ci{le}, \S 2.   

\medskip
In the  higher rank
case the curvature tensor is not, in general,  decomposable.
We introduce a notion of positivity which incorporates what is needed
to obtain  $L^2$-estimates-type results.

\medskip
Let   $\omega $ be  a hermitian metric
on $X$,  $\theta $ be a hermitian form on $T_X$ with continuous coefficients
and
$(X,E,\SS, h,h_s)$ be a s.h.m.; in particular, the curvature tensor
and the curvature form $\T{h}{E}$ are 
defined  a.e. (i.e. outside of $\SS$) and have measurable coefficients.

\begin{defi}
\label{ipo}
 {\rm ($\geq_t^{\mu}$; compare with  \ci{dem82}, \S 5.)}             
Let things be as above and $t$ be a positive integer. 
 
\noindent
We write: 
$$
\T{h}{E} \geq_t^{\mu} \theta \otimes {\rm Id}_{E_h}
$$ 
if
the following requirements are met.

\noindent
There exist a sequence 
of hermitian forms $\theta_s$ on $T_X\otimes E$  with continuous 
coefficients, 
 a sequence of continuous functions  $\lambda_s$ on $X$  
 and a continuous function $\lambda$ on $X$  subject to the following 
requirements:  

\noindent
{\rm (\ref{ipo}.1)}  $\forall x \in X$:
$|e_x|_{h_s}\leq  |e_x|_{h_{s+1}}$, $\forall s \in \nat$ and $\forall
e_x \in E_x$;

\noindent
{\rm (\ref{ipo}.2)}
$\theta_s\geq_t \theta \otimes {\rm Id}_{E_{h_s}}$;

\noindent
{\rm (\ref{ipo}.3)} $\T{h_s}{E}\geq_t \theta_s  - \lambda_s 
\omega \otimes {\rm Id}_{E_{h_s}}$;

\noindent
{\rm (\ref{ipo}.4)} $\theta_s  \to \T{h}{E}$ a.e. on $X$;

\noindent
{\rm (\ref{ipo}.5)} $\lambda_s \to 0$ a.e. on $X$;

\noindent
{\rm (\ref{ipo}.6)} 
$0\leq \lambda_s \leq \lambda$, $\forall s$.  
\end{defi}

Conditions  {\rm (\ref{ipo}.1)} and
{\rm (\ref{ipo}.6)} are needed   to apply Lebesgue's  
theorems
on monotonic and dominated convergence. 
In order to obtain $L^2$-estimates-type results, we also need 
the remaining  four conditions  to make precise the 
sought-for control of the curvature 
by the regularizing and approximating  metrics $h_s$.

\begin{rmk}
\label{refcoherent}
{\rm
 As an application of the $L^2$-estimates, we will see 
that  
 if  $\T{h}{E} \geq_N^{\mu} \theta \otimes {\rm Id}_{E_h}$,
for some continuous $\theta$, then
the sheaf  $\E(h)$ is coherent; see  Proposition \ref{coherent}. 
}
\end{rmk}

\begin{ex}
{\rm 
If $(E,h)$ is  a  hermitian bundle with $\T{h}{E} \geq_t \theta \otimes
{\rm Id}_{E_h}$, then  it is easy to exhibit $h$ as a s.h.m. such that 
$\T{h}{E} \geq_t^{\mu} \theta \otimes
{\rm Id}_{E_h}$; just   set
$h_s:=h$  $\forall s$, etc.
}
\end{ex}

\begin{ex}
\label{regappr}
{\rm 
Let $(E,h)$ be a vector bundle together with a {\em  continuous }
 s.h.m  metric. 
Under certain positivity conditions on the current
$id'd'' h^*$which is defined
on the total space of $E^*$ (see \ci{cm}, \S7.1) we
can exhibit $h$ as a s.h.m. with  positivity in the sense
of Definition \ref{ipo}. This is achieved in two steps.
In the first one
$h^*$ is regularized by
using Riemannian convolution coupled with the parallel transport
associated with an arbitrary hermitian metric on $E^*$
(see \ci{cm}, Lemme 7.2). In the second one the resulting metrics
are modified
so that they have the prescribed properties; this 
technical modification follows ideas
in \ci{dem82}, \S8. Details will appear elsewhere.  
}
\end{ex}

\begin{ex}
{\rm 
Let  $h=h_0e^{-2 \varphi}$ be  a singular metric  on a line bundle
$L$ 
with $T \geq \theta$ as currents
where $\theta $ is a continuous and real $(1,1)$-form. \ci{dem82}, 
Th\'eor\`eme 9.1,
exhibits   these data as   a s.h.m. $h$ with 
$\T{h}{L}\geq_1^{\mu}
\theta \otimes {\rm Id}_{L_{h}}$.

\noindent
Conversely, if we have  a s.h.m. $h$ with $\T{h}{L}\geq_1^{\mu} \theta
\otimes {\rm Id}_{L_h}$, then we have 
$\T{h_s}{L}\geq \theta_s - \lambda_s\omega$ and    
$T\geq  T_{ac} \geq \theta$.
}
\end{ex}

\begin{rmk}
{\rm 
The existence of a s.h.m. $h$ on 
a line bundle $L$ for which $\T{h}{L} \geq_1^{\mu} 0$    does not  imply that
$L$ is  nef.   
See {\rm \ci{d-p-s}}, {\rm Remark 1.6}.

\noindent
What is true is that if $L$ is nef, then $L$ will admit a metric
 $h=h_0e^{-2\varphi}$ with $h_0$  a hermitian metric on $L$ and $\varphi$
almost psh  such that
$\T{h}{L} \geq_1^{\mu} 0$.
 This can be seen by using
\ci{d-p-s},  Proposition 1.4, \ci{dem92}, Proposition 3.7 
and
 {\rm \ci{dem82}}, Th\'eor\`eme $9.1$.

\noindent
Similar remarks hold for big line bundles on projective manifolds
(cf. \ci{dem94},
 Proposition 6.6).
}
\end{rmk}

The following  lemma is elementary.

\begin{lm}
\label{fund}
Let  $(E,\SS_E, h,h_s)$ and
$(F,\SS_F, g,g_s)$  be s.h.m. on  two vector bundles
$E$ and $F$ over $X$, $\s_1$
and $\s_2$
be two real
 $(1,1)$-forms  with continuous coefficients
such that $\T{h}{E} \geq_{t_1}^{\mu} \s_1 \otimes {\rm Id}_{E_h}$
and
$\T{g}{F} \geq_{t_2}^{\mu} \s_2 \otimes {\rm Id}_{F_g}$.

\noindent
Then   $H:=h\otimes g$ on $E\otimes F$ can be seen as 
 a s.h.m. by setting $H_s:=h_s \otimes g_s$ and
$$
\T{H}{E\otimes F} \geq_{\min (t_1, t_2)}^{\mu}
 (\s_1 + \s_2)  \otimes {\rm Id}_{{(E\otimes F)}_H} \, .
$$
\end{lm}
Note that if the rank of $F$ is one, then 
$\min (t_1,t_2)= t_1$.

\medskip
 We  now prove that positivity is inherited by  quotient metrics.

\begin{lm}
\label{posquot}
Let $(X,E,\SS, h,h_s)$ be a s.h.m such that $\T{h}{E} \geq_t^{\mu} \theta 
\otimes
{\rm Id}_{E_h}$, $\phi:E\to Q$ be a surjection of vector bundles
with kernel $K$. Then $Q$ admits a s.h.m. $(Q, \SS' \subseteq \SS, q_s, q)$
 such that
 $\T{q}{Q}\geq_1^{\mu} \theta \otimes {\rm Id}_{Q_q}$.
\end{lm}

\noindent
{\em Proof.} 
Consider the dual exact sequence
$
0 \to Q^* \to E^* \to K^* \to 0.
$
 Each hermitian metric $h_s^*$ defines 
by restriction a hermitian metric $q_s^*$ on $Q^*$; analogously we
get $q^*:=h^*_{|Q^*}$. Clearly $(Q^*,\SS', q,q_s)$ is a s.h.m for an 
appropriate
$\SS'\subseteq \SS$. 
\noindent
For  every  $s$ we have that
 $\T{q_s^*}{Q^*}=\T{h_s^*}{E^*}_{|{Q^*}} + i\beta_s^* \wedge \beta_s$, 
where $\beta_s$  is a $(1,0)$-form with values in ${\rm Hom} (Q^*, 
K^*)$, 
${\cal C}^1$ coefficients and
$\beta^*$ is its adjoint. Moreover $i\beta_s^* \wedge \beta_s \leq_1 0$;
see \ci{dem82}, Lemme 6.6. The statement follows easily by dualizing 
again, 
which  has the effect of transposing and changing the signs.
\blacksquare

\medskip
More generally, a s.h.m. on $E$ ``with positivity" will induce s.h.m. 
``with positivity"
on $T^{\alpha}E $, $S^{\beta}E$ and $\wedge^{\gamma} E$. We leave the 
various 
formulations and elementary proofs to the 
reader.

\section{t-nef vector bundles}
\label{tnef}
\subsection{The definition of $t$-nefness}
\label{cn}

Let $X$ be a compact  manifold of dimension $n$, $\omega$ be  a hermitian 
metric
on $X$, 
$E$  be a  vector bundle of rank $r$ on $X$, $N:=\min (n,r)$
 and $L$ be a  line bundle
on $X$. Every tensor in $T_X \otimes E$ has rank $\rho \leq N$. 

\medskip
There are  notions of semi-positivity associated with every positive 
integer $t$. 
The standard one
is the following.

\begin{defi}
{\rm
(t-semi-positive vector bundle)}
We say that a vector bundle $E$ is  $t$-semi-positive, if 
$E$ admits a hermitian metric $h$ such  that $\T{h}{E} \geq_t 0$.
\end{defi}

Note that $E$ is $1$-semi-positive if and only if
it is  Griffiths-semi-positive, and that
$E$ is $N$-semi-positive if and only if it is Nakano-semi-positive. A 
similar remark
holds for strict inequalities.

\medskip
In algebraic geometry the most natural semi-positivity concept is {\em 
nefness}.
A differential-geometric characterization of this concept can be given as 
follows.

\begin{defi} 
\label{defneflb}
{\rm 
(Nef line bundle and nef vector bundle)}
We say that $L$ is nef if for every
$ \e >0$ there exists a hermitian metric $h_{\e}$ on $L$ such that
$\T{h_{\e}}{L}\geq -\e \omega$ as $(1,1)$-forms or, equivalently,
if $\T{h_{\e}}{L}\geq_1 -\e \omega \otimes {\rm Id}_{L_{h_{\e}}}$ as 
hermitian forms
on $T_X\otimes L$.

\noindent
We say that $E$ is  nef if the tautological line bundle
$\xi:=\odixl{{\Bbb P}(E)}{1}$ is nef.
\end{defi}
Note that the compactness of $X$implies that the  definitions given above
 are independent of 
the choice
of $\omega$. The same holds true for all  the other definitions given below
which involve a choice of $\omega$.

\medskip
If $X$ is projective, then  Definition \ref{defneflb} is equivalent 
to the usual one: {\em $L$ is nef if  $L\cdot C \geq 0$,
for every integral curve $C$ in $X$}; see  \ci{dem94}, Proposition 6.2.

\medskip
Unfortunately
   a nef line  bundle is not necessarily $1$-semi-positive ($\geq_1 0$). 
See {\rm \ci{d-p-s}, Example 1.7}, where an example  is given of a 
nef  rank two vector bundle $E$ on  an elliptic curve
such that 
the nef tautological line bundle $\xi$ on ${\Bbb P} (E)$ is not $\geq_1 0$.
Moreover, $E$  is not $\geq_1 0$ (otherwise $\xi$ would be
$\geq_1 0$); this shows that even on curves nefness and 
Griffiths-semi-positivity
do not coincide. 
Recall Theorem 1.12, \ci{d-p-s}, which 
states that nefness of a vector bundle $E$  can be characterized by the 
presence
of  
a system of hermitian metrics on all bundles
$S^{\alpha}E$ such that they are suitably semi-positive for
$\alpha \gg 0$.
I do not know if nefness can be characterized in terms of hermitian metrics
on the vector bundle itself.

The two facts above and the need to  express semi-positivity in terms of 
curvature
have 
motivated my   introducing   the notion of $t$-nefness.

\begin{defi}
{\rm
(t-nef vector bundle)}
We say that a vector bundle $E$ is  $t$-nef if for every $ \e>0$
there exists  a hermitian metric $h_{\e}$ on $E$ such  that 
$\T{h_{\e}}{E} \geq_t
-\e\omega \otimes {\rm Id}_{E_{h_\e}}$. 
\end{defi}

\smallskip
Every flat vector bundle  is $N$-nef. 

\noindent
If $E$ is $t$-semi-positive, then $E$ is $t$-nef. As pointed out above, 
the converse is 
not true in general; see  {\rm \ci{d-p-s}, Example 1.7}.

\noindent
If $E$ is $t$-semi-positive {\rm (}$t$-nef, respectively{\rm )}, then 
$E$ is
$t'$-semi-positive {\rm (}$t'$-nef, respectively{\rm )}, for every $ t'$ 
such that
  $1\leq t' \leq t$.

\noindent
By definition, a line bundle is nef if and only if it is $1$-nef.
A $1$-nef vector bundle is nef as we will see in {\rm Proposition 
\ref{list}}.
The converse is true on curves  as we will see in {\rm Theorem 
\ref{umemura}}.
 We do not know whether the converse is true or false when $\dim X \geq 2$.
This problem is
the analogue of Griffiths' question: does   ampleness  imply
Griffiths-positivity?

\noindent
We have checked that if $E$ is nef and 
is the tangent bundle of a compact complex surface or of a compact
K\"ahler threefold,  then $E$ is $1$-nef. This is done by
using the classification results contained in {\rm \ci{d-p-s}} and
{\rm Proposition \ref{list}}. According to conjectures in 
{\rm \ci{d-p-s}}, the same should be true
for compact K\"ahler manifolds of arbitrary dimension.

\noindent
({\bf From nefness to $1$-nefness})
On special manifolds, such as toric and abelian varieties, we have that if
 $E$  is nef, then $E \otimes \det E$ is  $1$-nef;
see {\rm \ci{cm}}, {\rm \S7.2.1}. By the following paragraph,
if $E$ is a rank $r$ nef vector bundle on such a variety, then 
$E \otimes (\det \, E)^{\otimes r+2}$ is $N$-nef.

\noindent
({\bf From $1$-nefness  to  $N$-nefness}) On any compact manifold,
if $E$ is $1$-nef, then $E\otimes \det E$ is $N$-nef. See
{\rm \ci{dem-sk}}.

\begin{ex}
\label{listnef}
{\rm 
({\bf Some $N$-nef vector bundles})
 The results mentioned above and the ones of sections \S\ref{bbpp} and  
\S\ref{nefoncurves}  give us the following list of examples.

\smallskip
\noindent
1) A nef vector bundle over a curve is $N$-nef. A nef line bundle
is $t$-nef for every $t$.

\noindent
2) A flat vector bundle is $N$-nef.

\noindent
3) 
If $X$ is a  special  manifold such as a toric or an 
abelian variety and $E$ is nef of rank $r$, 
then   ${E} \otimes \det {E}$ is $1$-nef and 
 ${E} \otimes (\det {E})^{\otimes r+2}$ is $N$-nef.

\noindent
4)
If $X$ is a K\"ahler manifold of dimension $n\leq 3$ with nef
tangent bundle $T_X$, then $T_X$ is $1$-nef and $K_X^{\otimes -1}\otimes T_X$
is $N$-nef.

\noindent
5)
Every Nakano-semipositive vector bundle
is $N$-nef (the converse is not true).
If $E$ is a 
 Griffiths-semipositive vector bundle, then  $E\otimes \det E$ is $N$-nef.

\noindent
6)
The extension of two $t$-nef vector bundles is $t$-nef. Positive tensor
representations of a $t$-nef vector bundle are $t$-nef. 
If $E_1$ and $E_2$ are $t$-nef, then  
$E_1\otimes E_2$ is $t$-nef. If $E$ is $t$-nef and $L$ is a nef line bundle,
then $E \otimes L$ is $t$-nef. 

\noindent
7) If $E$ is $1$-nef, then $E \otimes \det E$ is $N$-nef.

\noindent
8) If $f: X \to Y$ is a morphism and $E$ is a $t$-nef vector bundle on 
$Y$, then $f^* E$ is $t$-nef. If $t=1$ and, either $f$ is finite and 
surjective,
or $X$ is projective and $f$ has equidimensional fibers, then
the converse is true.

}
\end{ex}

At this point the following question is only natural.

\begin{??}
Is every nef vector bundle $1$-nef?
\end{??}

\subsection{Basic properties of $t$-nefness}
\label{bbpp}
Let us list and prove some basic properties of $t$-nef vector bundles.
We start with functorial ones.

\begin{pr}
\label{fctr}
Let $f:X \to Y$ be a holomorphic map, where $X$ and  $Y$  are  compact 
manifolds
and $E$ is a vector bundle on $Y$.

\medskip
{\rm (1)}
If $E$ if $t$-nef, then $f^*E$ is $t$-nef.

\medskip
{\rm (2)} Assume  that  $f$ is surjective 
and
 that the   rank of $E$ is one. 

\noindent
Then $f^*E$ is $1$-nef
{\rm (}=\,nef {\rm )} if and only if  $E$ is $1$-nef 
{\rm (}=\,nef {\rm )}.

\medskip
{\rm (3)} 
Assume that
$f$ is finite and surjective, that  $Y$ (and thus $X$)
is K\"ahler and let
 $E$ be of any rank.
Then  $f^*E$ is $1$-nef if and only if  $E$ is $1$-nef.  
\end{pr}

\noindent
{\em Proof.}
\medskip
\noindent
 (1). 
Let $\omega$ and $\omega'$ be two hermitian metrics on $X$ and $Y$, 
respectively. Let $A$ be a positive constant such that
$A \omega \geq f^* \omega'$. Fix $\e >0$ and let $\e':=\frac{\e}{A}$.
Let $h'$ be a hermitian metric on $E$ such that
$\T{h'}{E} \geq_t -\e' \omega' \otimes {\rm Id}_{E_{h'}}$. Endow $f^*E$ with
the pull-back metric $h:=f^*h'$. The claim follows from
the formula
$\T{h}{f^*E}= f^* \T{h'}{E}$.

\medskip
\noindent
{\rm (2).} See \ci{d-p-s}, Proposition 1.8.ii for the case of equidimensional
fibers and \ci{pm} for the general statement. 

\medskip
\noindent
(3)  It follows easily from 
 \ci{cm}, \S7.1: assign to $E$ the appropriate  trace metrics  and 
regularize.
\blacksquare

\begin{rmk}
{\rm 
As pointed out in {\rm Example \ref{regappr}}, the regularizing metrics 
in (3)
can be chosen
to satisfy  favorable conditions towards $L^2$-estimates.
}
\end{rmk}

\begin{??}
{\rm 
 Can we drop the assumption of finiteness from
 (3). A. J. Sommese
has pointed that the answer is 
positive when
$X$ is projective and the fibers are equidimensional: slice $X$ with 
sufficiently
ample general divisors to reduce to the case in which the morphism is 
finite. 
Is  (3)    true if we replace $1$-nef by
$t$-nef, with $t>1$? 
}
\end{??}

\begin{pr}
\label{list}
Let $X$, 
 $E$ and $r$ be as above. Then:

\medskip
{\rm (1)}
Let $E \to Q$ be a  surjection of vector bundles. If 
$E$ is $1$-nef,  then $Q$ is $1$-nef.

\medskip
{\rm (2)}
If $E$ is $1$-nef, then $E$ is nef.

\medskip
{\rm (3)} If $S^m E $ is $1$-nef, then $E$ is nef.

\medskip
{\rm (4)}
Let $0 \to K \to E \to Q\to 0$ be an exact sequence of vector bundles.
If  $K$ and $Q$ are  $t$-nef, then $E$ is $t$-nef.

\medskip
{\rm (5)}
Let $E=E_1 \oplus E_2$. The vector bundle $E$ is $t$-nef if and only if
$E_1$ and $E_2$ are $t$-nef.

\medskip
{\rm (6)} Let $F$ be another vector bundle.
Assume that  $E$  and $F$ are  $t$-nef and  $t'$-nef respectively;
then $E \otimes F$ is $\min (t, t' )$-nef.

\medskip
{\rm (7)}
 Assume that  $E$ is
$t$-nef.

\noindent
Then $S^m E$ and $\wedge^l E $ are $t$-nef
 for all $ m\geq 0$ and
for $ 0\leq l \leq r$.

\noindent
Moreover, $\Gamma^a E$ is $t$-nef, where
 $\Gamma^a E$ is  the irreducible tensor representation
of $Gl(E)$ of highest weight $a=(a_1, \ldots, a_r) \in {\zed}^r$,
with $a_1\geq \ldots a_r \geq 0$.

\medskip
{\rm (8)}
 Let $0 \to E \to E' \to \tau \to 0$ be an exact sequence with
 $E'$ a vector bundle on $X$  and $\tau$  a sheaf,
quotient of a $1$-nef vector bundle $E''$.
If $E$ 
is $1$-nef, then so is $E'$.

\medskip
{\rm (9)}
Let $0 \to K \to E \to Q\to 0$ be an exact sequence of vector bundles.
If
$E$ and $\det Q^*$  are  $1$-nef, then $K$ is $1$-nef.

\medskip
{\rm (10)}
Assume that $\det E$ is hermitian flat; the vector bundle
 $E$ is $t$-nef if and only 
if $E^*$ is $t$-nef.

\medskip
{\rm (11)} Let $E$ be $1$-nef and $s\in \Gamma (E^*)$. Then $s$ 
has no zeroes.
\end{pr}

\noindent
{\em Proof.} Fix, once and for all, $\omega$ a hermitian metric on $X$.

\medskip
\noindent
(1). Let $\e >0$ and  $h_{\e}$ be a hermitian metric on $E$
with $\T{h_{\e}}{E} \geq_1 -\e \omega \otimes {\rm Id}_{E_{h_{\e}}}$. Endow
$Q$ with the quotient metric $h'_{\e}$;  $Q$ can be seen as a
smooth 
sub-bundle of $E$  via the ${\cal C}^{\infty}$
orthogonal splitting of $E\to Q$ determined by $h_{\e}$, so that
${h_{\e}}_{|Q}=h'_{\e}$.
It is well-known 
(e.g. \ci{dem82}, Lemme 6.6)  that 
$\T{h'_{\e}}{Q} \geq_1 {{\T{h_{\e}}{E}}_{|Q}}_{h'_{\e}}$ 
and it is clear that 
${ 
  {\rm Id}_{ 
            {E_{ h_{\e} } }}}_{|Q}
={\rm Id}_{Q_{h'_{\e}}}$. The claim follows.

\medskip
\noindent
(2). Let $\pi : {\Bbb P}(E) \to X$ be the canonical projection. By virtue
of \ref{fctr}.1 we have that
 $\pi^* E$ is $1$-nef;   (1)
 and the canonical surjection
$\pi^*E \to  \odixl{{\Bbb P}(E)}{1}$ imply that this latter line bundle
is $1$-nef.

\medskip
\noindent
(3) $\pi^* S^m E$ is $1$-nef by \ref{fctr}.1, so that 
$\odixl{{\Bbb P}(E)}{m}$, being a quotient
of $\pi^* S^m E$, is $1$-nef, by (1). It follows 
that $\odixl{{\Bbb P}(E)}{1}$ is $1$-nef and thus nef.

\medskip
\noindent
(4). Fix a ${\cal C}^{\infty}$ vector bundle isomorphism $\Phi: 
E \to K \oplus Q$. Let  
and $\e >0$.
By assumption, there are metrics
$h_{K,\e}$ and $h_{Q,\e}$ such that
$\T{h_{K,\e}}{K} \geq_t -\frac{\e}{3} \omega \otimes {\rm 
Id}_{K_{h_{K,\e}}}$ and
$\T{h_{Q,\e}}{Q} \geq_t -\frac{\e}{3} \omega \otimes {\rm 
Id}_{Q_{h_{Q,\e}}}$.

\noindent
Fix an arbitrary positive real number $\rho >0$ and consider the
automorphism $\phi_{\rho}: Q \to Q$ defined by multiplication by the factor
$\rho^{-1}$.

\noindent
Let $\Phi_{\rho}:=({\rm Id}_K\oplus \phi_{\rho})\circ \Phi: E \to K 
\oplus Q$;
denote the first component of $\Phi_{\rho}$ by $\Phi_{K,\rho}$ and 
the second one by $\Phi_{Q,\rho}$

\noindent
Define a hermitian metric on $E$ by setting  $h_{\e, \rho}:=$
$\Phi_{K,\rho}^*h_{K, \e} \oplus \Phi_{Q,\rho}^*h_{Q, \e}$. Its 
associated 
Chern  connection has the form:

\[ D_{h_{\e, \rho}} =
\left(
\begin{array}{cc}
D_{h_{K,\e}}  &  -  \beta^*_{\rho} \\
\beta_{\rho}  & D_{h_{Q,\e}}
\end{array}
\right ),  \]

\noindent
where 
$\beta_{\rho}=\rho \beta_1$ 
is a 
$(1,0)$-form
with values
in 
${\rm Hom}(K,Q)$. By calculating $D^2$ we see that
$$
\T{h_{\e, \rho}}{E} \geq_t  -\frac{2}{3}\e \omega \otimes
{\rm Id}_{E_{h_{\e, \rho}}} + O(\rho)\omega \otimes {\rm Id}_{E_{h_{\e, 
\rho}}}.
$$
The claim follows by recalling that $X$ is compact
and by taking $\rho$ sufficiently small.

\medskip
\noindent
(5). The  ``if" part follows  from (4). 
The converse follows by observing
that if $E$ has a metric $h$, then each $E_i$ inherits a metric $h_i$ 
for which $D_{h_i}={D_{h}}_{|E_i}$. The same holds for the curvature 
tensors.

\medskip
\noindent
(6). The proof is immediate once one recalls the formula
for the curvature of the tensor product of two
hermitian metrics: $\T{h_1 \otimes h_2}{E_1\otimes E_2}=$
$\T{h_1}{E_1}\otimes {\rm Id}_{E_2}+$ ${\rm Id}_{E_1} \otimes 
\T{h_2}{E_2}$.

\medskip
\noindent
(7).
The tensor powers $T^n(E)$ are $t$-nef by virtue of
 (6). $S^n(E)$ and $\wedge^n(E)$
are both direct summands of $T^n(E)$ so that they are $t$-nef by
(5).
Recall that
$\Gamma^a E$ is a direct summand  of the vector bundle
$\otimes_{i=1}^r S^{a_i}(\wedge^i E)$ which is 
$t$-nef by what above, (5) and (6).

\medskip
\noindent
{\rm (8).}
The ``pull-back" construction 
gives the following
commutative diagram of coherent sheaves:
$$
\begin{array}{lllllllll}
 \hspace{1cm} 0 &
 \to  & 
E & 
\to   &
 E'  &
\to &
 \tau & 
\to &
 0  \\
\hspace{1cm} \, & \ & \uparrow {\rm Id}_E& \ & \uparrow q  & \  
& \uparrow & \ & \, \\
  \hspace{1cm} 0            & \to  & E & \to & 
E''' & \to  &  E''
&\to  & 0,              
\end{array}
$$

\noindent
where $q$ is surjective. Since $E$ and $E''$ are locally free and
$1$-nef, so is $E'''$ by $(4)$. Since $q$ is surjective,
it follows that $E'$ is $1$-nef by $(1)$.

\smallskip
\noindent
The proofs of (9) and (10) are the same as in 
the  nef case; the proof of (11) is in fact easier. 
The reader can consult
\ci{d-p-s}. 
\blacksquare

\begin{rmk}
\label{etl}
{\rm
It is easy to show, using 
$(6)$,  that if $E$ is $t$-nef and $L$ is a positive line bundle, then
$E\otimes L$ admits a hermitian 
metric $h$ with curvature $\Theta_h (E\otimes L) >_t 0$.
In particular, if $E$ is $N$-nef, then $E \otimes L$ is Nakano-positive.
A similar remark holds for the symbol $\geq_t^{\mu}$; see 
{\rm Lemma \ref{fundamental}}.
}
\end{rmk}

\begin{rmk}
{\rm
As far as $(1)$ above is concerned, it is not true that if
$E$ is  $t$-nef, then  $Q$ is  $t$-nef. In fact, consider the 
canonical surjection $\odix{\pn{2}}^3 \to T_{\pn{2}}(-1)$: 
$\odix{\pn{2}}^3$ is 
$2$-nef, but if $T_{\pn{2}}(-1)$ were $2$-nef, then $T_{\pn{2}}=$ 
$T_{\pn{2}}(-1)
\otimes \odixl{\pn{2}}{1}$ would be $>_2 0$,  i.e.
 Nakano-positive and this is 
a contradiction. This example also shows that $1$-nefness
is strictly weaker than $2$-nefness.
We do not know whether $(8)$ is false when we replace $1$ by $t$. 
}
\end{rmk}

\subsection{Nefness and $t$-nefness on curves}
\label{nefoncurves}
It is an outstanding problem in Hermitian differential geometry
to determine whether an ample vector bundle is Griffiths-positive.
In
\ci{ume}, Umemura proves that on curves 
 ampleness and Griffiths-positivity 
coincide. As is was pointed out to me by N.M.
Kumar, the part of the argument that needs a result
analogue to Proposition
\ref{list}.8 
is omitted in \ci{ume}.
 
\medskip
We now prove  that on curves  nefness
 and $1$-nefness coincide: 
 the algebraic notion  of nefness can be characterized
in differential-geometric terms. Recall that Example 1.7, \ci{d-p-s},
implies that even over a curve,  a nef vector bundle is not necessarily
Griffiths-semi-positive. 

\begin{tm}
\label{umemura}
Let $X$ be a nonsingular projective curve
and  $E$ be a vector bundle of rank $r$ on $X$.
The following are equivalent.

\smallskip
{\rm ($i$)} $E$ is $1$-nef;

\smallskip
{\rm ($ii$)} $E$ is nef;

\smallskip
{\rm ($iii$)} every quotient bundle of $E$, and in particular
$E$,  has non-negative degree. 
\end{tm}

\noindent
{\em Proof.}
($i$)$ \Rightarrow$($ii$). This is Proposition \ref{list}.2.

\smallskip
($ii$)$ \Rightarrow$($iii$). In fact they are equivalent by \ci{c-p},
Proposition 1.2.7.

\medskip
($iii$)$ \Rightarrow$($i$). We divide the proof in three cases, according
 to whether $g=0$, $g=1$ or $g\geq 2$. Let $d$ be the degree of $E$.
By assumption $d\geq 0$.

\medskip
If the genus $g(X)=0$, then $E$ splits into a direct 
sum of  
line bundles and the statement follows  easily.

\medskip
Let $g(X)=1$. It is enough to consider
the case when $E$ is indecomposable.
Let us first assume that $d\geq r$.  By
 \ci{at}, Lemma 11,  $E$  admits a maximal splitting
$(L_1,\ldots , L_r)$ with $L_i$ ample line bundles on $X$. It follows that
$E$ could then be constructed inductively from (ample =) positive line 
bundles
by means of extensions. A repeated use of Proposition
\ref{list}.4 would allow us to conclude.
We may  thus assume, without loss of generality,
that $0\leq d <r$. 
If $r=1$, then $E$ is either ample
or hermitian flat; in both cases we are done. We now proceed by induction on
the rank of $E$. Assume that we have proved our contention
for every vector bundle of rank strictly less than $r$. By \ci{at}, Lemma
15 and Theorem 5, $E$ sits in the middle of an exact sequence:
$$
0 \to A \to E \to B \to 0,
$$ 
where $B$, being a quotient of $E$,
enjoys property ({\em iii}) and
$A$ is either  a trivial vector bundle (if $d>0$) or a
hermitian flat line bundle
(if $d=0$). In any case $A$ is clearly
$1$-nef and $B$ is $1$-nef by the induction hypothesis.
 We can apply
\ref{list}.4  and  conclude that $E$ is $1$-nef. 
This proves the case $g(X)=1$.

\medskip
We now  assume that $g(X)\geq 2$. The proof will be by induction
 on $r$.
If $r=1$, then  we are done since $\deg{E}\geq 0$ implies that either $E$ 
is ample
or it is hermitian flat. Assume that we have proved our assertion
 for every vector bundle of rank strictly less than $r$.

\smallskip
\noindent
There are two cases. 

\noindent
In the first one we
suppose that $E$ contains a non-trivial  vector  sub-bundle 
$K$ which is $1$-nef. Consider
the exact sequence of coherent shaves:
$$
0 \to K \to E \to Q:=E/K \to 0.
$$
There are two sub-cases. In the first one we assume that
 $Q$ is locally free.
By assumption every quotient vector bundle of $Q$, being in turn
a quotient bundle of $E$,   has positive degree. The induction hypothesis
forces $Q$ to be  $1$-nef. Proposition \ref{list}.4 allows us to conclude 
that
$E$ is $1$-nef as well.

\noindent
In the second sub-case $Q\simeq F  \oplus \tau$, with $F$ locally free 
and $\tau$ 
has zero-dimensional support; in particular   there is a surjection
$\odix{X}^m \to \tau$. 
If $K'$ is the kernel of the surjection $E\to F$, then 
we have the exact sequence
$$
0 \to K \to K' \to \tau \to 0,
$$
so that, by \ref{list}.$8$, $K'$ is $1$-nef and we are reduced to the 
first sub-case. 

\noindent 
In the second case we are allowed to assume that $E$ does not  contain 
properly
  any
non-trivial vector bundle $K$ which is $1$-nef. 

\smallskip
{CLAIM.}
$E$ is  stable. 

\noindent 
To prove this we start a new proof by induction.
Let $K$ be a vector bundle contained in $E$, neither trivial nor
equal to  $E$. Since (iii) implies that $\deg{(E)}\geq 0$, to prove that
$E$ is stable it is enough to show that $\deg {K} <0$.
Seeking a contradiction, let us assume that $\deg {K} \geq 0$.
Let $s$ be the rank of $K$. If $s=1$, since  $K$ is not  $1$-nef by the 
working assumption of this second case, we see that $\deg{K}<0$ (otherwise
$K$ would be either ample or hermitian flat) and we have reached
a contradiction if $s=1$.  Assume that, for every
non-trivial 
$K''\subseteq E $ with rank strictly less than $s$, $\deg{K''} <0$. 
Each sub-bundle
of $K$ has, by this second inductive hypothesis, negative degree.
Since we are assuming that  $\deg {K} \geq 0$, it  follows that every 
quotient
of $K$, including $K$ itself has  non-negative degree, so that, by the 
first 
induction hypothesis,
$K$ is $1$-nef and we have reached a contradiction for every $s$: 
the degree of $K$ must 
be negative,
 $E$ is stable and the claim is proved.

\smallskip
\noindent
 Since $\det E$
has non-negative degree by assumption,
$\det E$ is $1$-nef. Since $E$ is stable, for any hermitian metric $h$ on 
$\det E$  of curvature $\T{h}{\det E}$,
 a standard calculation
(see \ci{ume}, Lemma 2.3) yields  a hermitian metric $H$ on $E$
of curvature $\frac{1}{r}\T{h}{\det E} \otimes {\rm Id}_E$. This proves that
$E$ is $1$-nef also in the second case.
\blacksquare

\subsection{A differential-geometric characterization of ampleness for 
vector bundles}
\label{dg}
We  now prove a characterization of ampleness by means of curvature
 properties  which is a simple consequence of \ci{d-p-s}, Theorem 1.12.

\begin{pr}
\label{1.12}
Let $X$ be a compact manifold equipped with a hermitian metric $\omega$
and $E$ be a vector bundle on $X$.
Then $E$ is ample if and only if there exists a sequence
of hermitian metrics  
$h_m$ on $S^m E$ such that 

\smallskip
{\rm (i)} the sequence of metrics on  
$\odixl{{\Bbb P}(E)}{1}$ 
induced by the surjective 
morphisms 
$$\pi^* S^m E \to \odixl{{\Bbb P}(E)}{m}$$
 converges uniformly
to a hermitian  metric $h$ of positive
curvature on  $ \odixl{{\Bbb P}(E)}{1}$ and

\smallskip
{\rm (ii)}
there exist $\eta >0$ and $m_0 \in \nat$
such that
 $\forall m \geq m_0$:
$$
\label{F&A}
\T{h_m}{S^m(E)} \geq_1 m\eta \omega \otimes {\rm Id}_{{S^m E}_{h_m}}.
$$
\end{pr}

\begin{rmk}
\label{known}
{\rm
If $E$ is ample, then
the fact that  some metrics $h_m$ 
with property  (ii) exist  for all $ m \gg 0$ is a 
well known  consequence of
{\rm \ci{griff}}, theorems F and A. The point  made by the statement above
 is that the metrics
$h_m$ are constructed on {\em  all} symmetric powers $S^m(E)$,
 and that they are all  built starting from
a suitable metric on $\odixl{{\Bbb P}(E)}{1}$; see {\rm  \ci{dem92}, 
Theorem 4.1}. 
}
\end{rmk}

\noindent
{\em Proof.} The proof of the implication ``$\Leftarrow$" follows 
easily from (i): $ \odixl{{\Bbb P}(E)}{1}$ is positive by the existence of
$h$ so that 
$E$ is ample.

\noindent
For the reverse implication ``$\Rightarrow$" we argue as follows.
Fix a hermitian metric $\omega'$ on ${\Bbb P}(E)$. The ampleness of
$E$ implies the ampleness of
$\odixl{{\Bbb P}(E)}{1}$, which then  admits a hermitian metric $h$
of positive curvature;  the compactness of $X$
 ensures us that there exist $\alpha >0$  and  $A>0$ such that
$$
\T{h}{\odixl{{\Bbb P}(E)}{1}} \geq \alpha 
\omega' \geq \alpha A \pi^* \omega.
$$
Define $\eta:=\frac{2}{3} \alpha A$. We are now in the position of using
\ci{dem92}, Theorem 4.1  with $v:=\frac{3}{2}\eta \omega$ and 
$\e:=\frac{1}{2}\eta$.
\blacksquare

\section{Vanishing theorems}
\label{vt}
In this section
 we link the positivity of $h$  to the vanishing of the cohomology
of $K_X \otimes \E (h)$.

\subsection{The basic $L^2$-estimate and vanishing theorem, and the coherence
of $\E (h)$}
Following Demailly, \ci{dem82}, \S 5, we say that  a s.h.m. 
$(X,E,\SS,h,h_s)$ is
{\em $t$-approximable} if $\T{h}{E} \geq_t^{\mu} 0$ (cf. Definition
\ref{ipo}).
We denote the space of $(p,q)$-forms with values in $E$ and coefficients
which are locally square-integrable
 by $L_{p,q}^2(X,E,{\rm loc})$. As usual, $n=\dim{X}$, $r$
is the rank of $E$ and $N:={\rm min} (n,r)$.

\begin{pr}
\label{l2}
{\rm (See \ci{dem82}, Th\'eor\`eme 5.1)}
Let $(X,\omega)$ be
  K\"ahler, where either 
 $\omega$  is  complete or $X$ is weakly pseudoconvex.
Assume that  $(E,h)$ is a s.h.m. with the property that
 $\T{h}{E}\geq_{n-q+1}^{\mu}
\e \omega \otimes {\rm Id}_{E_h}$, where $\e$ is a non-negative and 
continuous function on $X$ and $q>0$ is a positive integer.

\noindent
Let  $g\in L^2_{n,q}(X,E,{\rm loc})$ be such that 
$$
d''g=0,        \qquad \qquad   
\int_X{  |g|^2_h \, dV_{\omega}} < + \infty 
\qquad  and \qquad 
 \int_X{ \frac{1}{\e} |g|^2_h \, dV_{\omega}} < + \infty .
$$
\noindent
Then there exists $f \in L^2_{n,q-1}(X,E,{\rm loc})$ such that
$$d''f=g \qquad and  \qquad
 \int_X  |f|^2_h dV \leq
\frac{1}{q} \int_X \frac{1}{\e} |g|^2_h \,  dV_{\omega}.
$$
\end{pr}

\noindent
{\em Sketch of  proof.} Th\'eor\`eme 5.1 states something
slightly different but it is immediate to recover the statement
of the proposition. We merely point out, for the reader's convenience,
 the minor changes to be implemented to obtain the above statement. The 
notation
is from \ci{dem82}.

\noindent
The assumption $\T{h}{E}\geq_{n-q+1}^{\mu}
\e \omega \otimes {\rm Id}_{E_h}$ has two consequences. The former is that
$h$ is $n-q+1$-approximable. The latter is that, by virtue of
\ci{dem82},
Lemme 3.2 (3.4):
$$
|g|^2_{\T{h}{E}} \leq  \frac{1}{q\e} |g|^2_h \quad a.e.
$$
We can apply the aforementioned theorem and conclude.
\blacksquare

\medskip
The following generalizes Nadel Vanishing Theorem.
It is an easy consequence of the proposition above.

\medskip
\begin{tm}
\label{vanish}
Let $(X, \omega)$ be K\"ahler with $X$ weakly pseudoconvex.
Assume that  $(E,h)$  is  a s.h.m. such that 
$\T{h}{E} \geq_{N}^{\mu} \epsilon \omega \otimes {\rm Id}_{E_h}$ 
for some positive and continuous function
$\epsilon$.
Then,
$H^q(X,K_X\otimes \E (h))=0$,  $\forall
q>0$.
\end{tm}

\noindent
{\em Proof.}
 The complex $({\frak L}^{\bullet}, d'')$ of \S\ref{mi} is exact by 
Proposition
\ref{l2} applied to small
balls. 
This complex is therefore an acyclic resolution of
$K_X \otimes \E(h)$ whose  cohomology is isomorphic to the cohomology
of the complex of global sections of  $({\frak L}^{\bullet}, d'')$.
This latter cohomology is trivial for every positive value of $q$ 
by  Proposition
\ref{l2} (modify the metric as in \ci{dem94}, Proposition 5.11). 
\blacksquare
\medskip

\medskip
We now prove that 
   if  $h$ is suitably positive,  then $\E (h)$ is  coherent.
The line bundle case is due to Nadel.

\begin{pr}
\label{coherent}
Let $X$ be a complex manifold,
$(X,E,\SS,h,h_s)$ be a s.h.m.,
 and $\theta$ be a  continuous real $(1,1)$-form
on $X$ such that $\T{h}{E} \geq_N^{\mu} \theta \otimes {\rm Id}_{E_h}$.
 Then $\E (h)$ is coherent.
\end{pr}

\noindent
{\em Proof.} We make the necessary changes from the line bundle
case (cf. \ci{dem94}, Proposition 5.7). 

\noindent 
Note that the condition 
$\T{h}{E} \geq_N^{\mu} \theta \otimes {\rm Id}_{E_h}$ 
implies that  $h\geq  h_1$ a.e.

\noindent
The statement being local, we may assume that 
$X$ is a ball centered
about the origin in $\comp^n$
with holomorphic coordinates $(z)$,
that $E$ is trivial
and that $\theta$ has bounded coefficients.
 Let $\omega$ be the $(1,1)$-form associated with
 the euclidean metric on $X$.
 Let $\frak S$ be the vector space of holomorphic sections 
$f$ of $E$ such that $\int_X{|f|^2_h \, d\lambda}  < \infty$, where 
$d\lambda$  is  the Lebesgue measure on $\comp^n$. Consider
 the natural evaluation map $ev:{\frak S} \otimes_{\comp} \odix{X} \to E$.
The sheaf ${\frak E}:= Im(ev)$  is coherent by Noether
Lemma (cf. \ci{gr-re}, page 111) and it is contained in $\E (h)$.

\noindent
We want  to prove that $\E (h)_x = {\frak E}_x$ for all $x \in X$.  
In view of Nakayama's Lemma,
\ci{at-mac}, Corollary 2.7, it is enough to show that
$ {\frak E}_x + {{\frak m}^{\gamma}_x} \cdot \E (h)_x = \E (h)_x$ for 
some $\gamma
\geq 1$.

\medskip
\noindent
STEP I. Assume that we could prove that:

\smallskip
\noindent
$(\bullet) \qquad \qquad \qquad \qquad $ 
$ {\frak E}_x + \E (h)_x \cap {\frak m}^l_x  \cdot E_x =\E (h)_x$ for 
every positive
integer $l$. 

\smallskip
\noindent 
By the Artin-Rees Lemma, \ci{at-mac}, Corollary 10.10, there would 
be  a positive
integer $k=k(x)$ such that
$$
\E(h)_x = {\frak E}_x + \E(h)_x\cap {\frak m}_x^l \cdot E_x \subseteq
{\frak E}_x + {\frak m}_x^{l-k}\cdot \E(h)_x 
\subseteq {\frak E}_x + {\frak m}_x \cdot \E(h)_x
\subseteq \E(h)_x
$$
 for all $l\geq k$. All symbols ``$\subseteq$" could be replaced by
equalities and we could conclude that 
${\frak E}_x = \E (h)_x$ by Nakayama's Lemma as above.

\medskip
\noindent
STEP II. We now prove $(\bullet)$.

\noindent
Let $f$ be a germ in $\E (h)_x$ and $\s$ be a smooth cut-off function 
such that  is identically $1$ around $x$ and that 
has compact 
 support small enough so that $\s f$ is smooth  on $X$.

\noindent
For every positive integer $l$
define  a strictly psh function $\varphi_{l}:= (n+l)\ln |z-x| +C|z|^2$
where $C$ is a positive constant chosen so that
$2id'd'' (C|z|^2)  +  \theta \geq \e \omega$, for some positive constant 
$\e$.

\noindent
Define a metric on $E$ by setting
$H_l:=$
 $h e^{-2\varphi_{l}}$.
Since both $\ln |z-x|$  and $|z|^2$ are  psh, we can apply
the results of \ci{dem82}, \S9 to $\varphi$ and deduce, with the aid of Lemma
\ref{fund}, that  $H_{l}$ is a s.h.m. on $E$ with
 $\T{H_{l}}{E}=\T{h}{E} 
+ 2id'd''\varphi_{l} \otimes {\rm Id}_E$  
and such  that   $\T{H_{l}}{E} \geq_N^{\mu} \e \omega \otimes
{\rm Id}_{ E_{H_{l}} }$. 

\noindent
Consider the smooth  $(0,1)$- form $g:=d''(\s f)$ which has compact support
and is
identically zero around $x$. The function $|z-x|^{-2n -2l}$
is continuous outside $x$. It follows that:
$$
\int_X{|g|_{H_{l}}^2\, d \lambda}= 
\int_X{ 
       |g|_h^2 \,  
                 |z-x|^{-2n- 2l}\, e^{-2C|z|^2} \,
                                                  d \lambda < \infty.
} 
$$
We solve, for every index $l$, the equation $d''u=g$ with $L^2$-estimates
 relative to $H_{l}$ using
Proposition \ref{l2}. We obtain  a  set of solutions $u_{l}$ such that
$$
\int_X{|u_{l}|^2_{H_{l}}} \, d\lambda=
\int_X{|u_{l}|^2_{h} \,  |z-x|^{-2n-2l} \,  e^{-2C |z|^2} \, d\lambda
 < \infty}.
$$
Since the factor $e^{-2C |z|^2}$ does not affect integrability
we get that
$$
\int_X{|u_{l}|^2_{h} \,  |z-x|^{-2n-2l}  \, d\lambda
 < \infty}.
$$
Since
$d'' (\s f - u_{l}) =0$ and $h \leq H_{l}$,  we see  that 
$ \s f- u_{l}=: F_{l} \in {\frak E}$ (cf. \ci{g-h} page 380). 
The germ $u_{l,x}=f  - F_{l,x}$ is holomorphic.
  Since $h\geq h_1$ and $h_1$ is continuous, there is a positive constant
$B$ such that:
$$
 \int_X{ B \,  |u_{l}|^2  \,  |z-x|^{-2n-2l} } \, d\lambda 
\leq 
\int_X{|u_{l}|^2_h  \, |z-x|^{-2n-2l}  } \, d\lambda
 < \infty.
$$ 
Let $u_{l}^{\{j \}}$ be the $j$-th coordinate  function of $u_{l}$,
$j=1, \ldots ,r$.
By a use of Parseval's formula (cf. \ci{dem94}, 5.6.b) we see
that $u_{l}^{\{j \}} \in {\frak m}^l_x$ for every index $j$.
It follows that $(\bullet)$  
holds and we are done.
\blacksquare

\subsection{$t$-nefness and vanishing}
We now show how to use Theorem \ref{vanish} to infer
the vanishing of cohomology in the case
of a $N$-nef vector bundle twisted by 
 a line bundle which  can be endowed with a positive s.h.m.

\smallskip 

The following is an elementary consequence
of Lemma \ref{fund}:

\begin{lm}
\label{fundamental}
Let $E$ be a $t$-nef vector bundle on a compact manifold $X$,
$\omega$ be a hermitian metric on $X$, 
$\theta$ be a real $(1,1)$-form with continuous coefficients and
$(F,g,g_s)$ be a vector bundle endowed with a s.h.m. such that
$\T{g}{F} \geq_t^{\mu} \theta \otimes {\rm Id}_{F_g}$.

\noindent
Then for every constant $\eta >0$ there is a  s.h.m. $H_{\eta}$ on 
$E\otimes F$
for which:
$$
\T{H_{\eta}}{E\otimes F} \geq_t^{\mu} (\theta - \eta \omega) \otimes
{\rm Id}_{(E \otimes F)_{H_{\eta}}}.
$$
Moreover, if $F$ is a line bundle and
$\T{g}{F} \geq \theta$ as $(1,1)$-forms, then the same conclusion holds.
\end{lm}

\begin{lm}
\label{ELco}
Let
 $(F,h_F)$ be a hermitian vector bundle  on  a manifold
$X$ and $(L,h_L)$
be a line bundle  on $X$ endowed with  a singular metric
 $h_L$ as in {\rm \S\ref{exlbshm}}.
Consider the vector bundle $E:=F\otimes L$ endowed with the
measurable metric
$h:=h_F\otimes h_L$.

\noindent
Then $ \E (h) =\id (h_L) \otimes E$ and
$\E (h)$ is coherent. 
\end{lm}

\noindent
{\em Proof.} The statement  $\E (h) = \id (h_L)\otimes E$ is local on $X$
 so that 
 we may assume that $X$ is a ball in $\comp^n$,  
 that $F$ and $L$ are trivial, that  $h_L=e^{-2\varphi}$
with $\varphi$ almost psh and that $h_F$ has bounded coefficients.

\noindent
Let us first prove that $ \id (h_L)_x \otimes E_x \subseteq \E(h)_x$.
Let $e_x\in  E_x$ and  $f_x \in \id (h_L)_x$.
Since $h_F$ is continuous,
we have that  $|f_x e_x|^2_h= |f_x|^2 |e_x|^2_{h_F} e^{-2\varphi}$ 
is locally integrable.

\noindent
Let us prove the reverse inclusion  $ \E (h)_x
\subseteq \id (h_L)_x\otimes E_x  $ for every $x$ in $ X$.
 There exists a constant $\tilde\epsilon >0$ 
 such that $h_F \geq \tilde\epsilon \Delta$, where
$\Delta$ is the standard euclidean metric on the fibers
of $E$. Fix $x \in X$.
Assume that $\E (h)_x \ni e_x=<f_1,\ldots, f_r>$.
Then $|e_x|^2_h=|e_x|^2_{h_F} e^{-2\phi} \geq
\tilde
\epsilon \sum |f_i|^2 e^{-2\varphi}$. As the  
left hand side of the inequality is integrable around $x$, so
is each summand on the right. This proves the reverse inclusion.
To conclude  recall that $\id (h_L)$ is coherent (or apply Proposition
\ref{coherent}).
\blacksquare

\medskip
The following result
is the key to
the proofs of 
the effective statements to be found in    \S\ref{effective}.
See Ex.  \ref{listnef} for examples of $N$-nef vector bundles.

\begin{tm}
\label{myvan}
Let $(X,\omega)$ be as in 
{\rm Theorem \ref{vanish}}, and 
  $(F,h_F)$, $(L,h_L)$ and $(E,h)$ be as in {\rm Lemma
\ref{ELco}}.

\noindent
If $\T{h}{E}\geq_N^{\mu} 
\e \omega \otimes
{\rm Id}_E$ for some positive and continuous function $\e$, then
$$
H^q(X, K_X  \otimes F \otimes L \otimes  \id (h_L) )=
H^q(X, K_X \otimes \E (h) ) =0, \qquad \forall \, q>0.
$$
Moreover, if $X$ is compact,
 $F$ is $N$-nef
and $(L,h)$ is such that  
 $\T{h}L\geq \e \omega$,  for some positive constant $\e$,
then the same conclusion holds.
\end{tm}

\noindent
{\em Proof.} By
Lemma \ref{ELco}, we have that $\E(h) =\id (h_L) \otimes E$. 
We conclude in view of  Theorem \ref{vanish}.

\noindent
The case of $X$ compact is a special case after Lemma \ref{fundamental}. 
\blacksquare

\medskip
The following is not needed in the sequel. We include it since it 
is a generalization
of Kawamata-Vieheweg Vanishing Theorem (K-V) and it can be proved
along the lines of \ci{dem94}, 6.12 by 
using Theorem \ref{myvan} instead of 
Nadel Vanishing Theorem.  The ``$1$-nef" case 
follows easily from
K-V Theorem  and the Leray spectral sequence
by looking at the projectivization of
$E$. The statement in the ``$N$-nef" case seems new for $0<q< {\rm rank} 
\, E$
and the vanishing in the complementary range  
follows from   K-V and Le-Potier spectral sequence.
 See
\ci{dem94} for the particular language employed in the statement below.

\begin{tm}
\label{kv}
Let $(X,E,F)$ be the datum of: $X$ a projective manifold, $E$ a 
$N$-nef
vector bundle on $X$, $F$ a line bundle on $X$ such that some 
positive multiple
$mF$ can be written as $mF=L+D$, where $L$ is a nef line bundle and $D$
is an  effective divisor.
Then
$$
H^q(X, K_X \otimes E \otimes F \otimes \id(\frac{1}{m}D))=0 \qquad 
for\, \, q
> \dim X - \nu (L),
$$
where $\nu(L)$ is the numerical dimension of $L$ and $\id(\frac{1}{m}D)$
is the multiplier ideal of the singular local weights
associated with the  $m-$roots of the absolute values of local  
equations for $D$. 

\noindent
As a special case, we have that if $F$ is a nef line bundle, then
$$
H^q(X, K_X \otimes E \otimes F )=0 \qquad for \, \, q
> \dim X - \nu (F);
$$
in particular, if $F$ is  nef and big, then 
$$
H^q(X, K_X \otimes E \otimes F )=0 \qquad for \, \, q
> 0.
$$
\end{tm}

\section{Effective results}
\label{effective}
\subsection{Special s.h.m. on line bundles after Anghern-Siu, Demailly, Siu
and Tsuji}
\label{demsiu}
The following proposition is  at the heart
of the effective base-point-freeness,  point-separation and 
jet-separation results
in \ci{an-siu}, \ci{tsu},
 \ci{siu-va}, \ci{siu94b} and \ci{dem96};
it  provides us with the necessary
s.h.m.  which we   transplant to the vector bundle case
and use in connection with Theorem \ref{myvan}.

\medskip
First we need to fix some notation. 

\bigskip
Let $F$ be a rank $r$ vector bundle on a complex manifold  $X$ and $p$ be any
 positive integer.
We say that {\em the global sections of $F$ generate simultaneous 
jets of order
$s_1,\ldots, s_p \in {\Bbb N}$ at arbitrary $p$ distinct 
points of $X$}
if the natural maps
$$
H^0(X,F) \to \bigoplus_{\i=1}^p  {\cal O}(F)_{x_{i}} 
\otimes \odix{X}/{\frak m}^{s_{i}+1}_{x_{i}}
$$
are surjective for every choice of $p$ distinct points  $x_1, \ldots , x_p$
in $X$.

\noindent
We say that 
{\em 
the global  sections of $F$ separate arbitrary $p$  
distinct points
 of  $X$
} 
if the above holds with
all $s_{i}=0$.

\smallskip
Assume that $X$ is compact. Let $V:= H^0(X,F)$
and $h^0:=h^0(X,F):= \dim_{\comp} H^0(X,F)$.
Consider $G:= G(r,h^0)$ the Grassmannian of $r$-dimensional quotients
of $V$, $\QQ$ the universal quotient bundle of $G$ 
and $\q$ the determinant of $\QQ$.

\noindent
As soon as $F$ is generated by its global sections
(which corresponds to the above conditions being met
for $p=1$ and $s_1=0$), we get a morphism
$f:X \to G$ assigning to each $x \in X$ the quotient $F_x \otimes k(x)$
and such that $F \simeq f^* \QQ$. The Pl\"ucker embedding defined by
$\q$ gives a closed embedding into the appropriate projective space
$\iota:  G \to {\Bbb P}$. We obtain a closed embedding 
$\hat f := \iota \circ f: X \to {\Bbb P}$.
It is clear that:

\noindent
- $V$ separates arbitrary $2$ points of $X$ iff  $f$ is 
bijective birational onto its image;

\noindent
- If $V$ separates arbitrary pairs of points of $X$ and generates
jets of order $1$
at an arbitrary point of $X$, then  $f$ is a closed embedding.

\smallskip
Given $n$, $p$ and $\{ s_1, \ldots , s_p \}$ as above let us
define the following integers:
$$
m_1 (n,p) : = \frac{1}{2}(n^2 +2pn -n +2 ),
$$
$$
m_2(n,p;s_1, \ldots , s_p)= 2n \sum_{i=1}^p B(3n + 2s_i -3, n) + 2pn +1,
$$
where $B(a,b)$ denotes the usual binomial coefficient,
$$
m_3(n,p;s_1, \ldots , s_p) = (pn  + \sum_{i=1}^p s_i) \, m_1(n,1)
$$
and
$$
m_4 (n)= (n+1) \, m_1(n,1).
$$
\begin{pr} 
\label{heart} 
Let $X$ be a projective manifold of dimension $n$ and $L$ be
an ample 
line bundle on $X$. Fix a K\"ahler form $\omega$ on $X$.

\medskip
\noindent
{\rm (\ref{heart}.1) (Cf. \ci{an-siu} and \ci{tsu}.)} 
Let  $p$ be a positive integer. Assume that 
$m
\geq m_1(n,p)$.

\noindent
 Then
for any set of $p$ distinct points $\{x_1, \ldots x_p\}$ of $X$,
there exists a nonempty 
subset $J_0\subset \{1,\ldots , p\}$ with the following property:

\noindent
there exist $\e>0$,  a s.h.m. $h$ for
 $mL$ with  
 $\T{h}{mL} \geq^{\mu}_1 \epsilon \omega \otimes {\rm Id}_{L_h}$
and with the property that 
the multiplier ideal $\id (h)$ of $h$ is such that the closed subscheme
 given by $\id (h)$ has the points $x_{i}$ as isolated points $\forall
i \in J_0$ and contains all the  points $\{ x_{i}\}$. 

\medskip
\noindent
{\rm (\ref{heart}.2) (Cf. \ci{siu94b}; see also  \ci{dem96}.)} 
Fix a positive integer $p$ and
 a sequence of non-negative integers $\{s_1, \ldots , s_p\}$. 
Assume that 
$
m\geq m_2(n,p; s_1, \ldots , s_p).
$

\noindent 
Then for any set of $p$ distinct points
 $\{x_1, \ldots
x_p\}$ of $X$
there exist $\e>0$,  a s.h.m. $h$ for
 $K_X+mL$ with  
$\T{h}{K_X + mL} \geq^{\mu}_1 \epsilon \omega  \otimes {\rm Id}_{L_h}$
and with the property that 
the multiplier ideal
$\id (h)$ satisfies   $\id(h)_{x_{i
}} \subseteq {\frak m}_{x_i}^{s_i+1}$, for every
$1\leq i \leq p$, and
is such that the closed subscheme
 given by $\id (h)$ has all the points $x_{i}$ as isolated points.
\end{pr}

The easy lemma that follows is probably well-known and
  makes precise a well-understood principle: {\em
it is easy to go from global generation to the generation of higher jets}.
Though  the presence of the nef line bundle  $M$ is redundant
in the statement,  we use it  because
of the application of this lemma
 to the case of higher rank.

\begin{lm}
\label{freetojet}
{\rm (From freeness to the generation of jets)}
Let $X$, $n$, $p$ and $\{ s_1, \ldots , s_p \}$ be as above, 
 $F$, $A$ and $M$  be   line bundles on $X$ such that
$F$ is  ample and generated by its global sections, $A$ is ample
and $M$ is  nef.
Then the global sections of $K_X + (pn + \sum_{i=1}^p s_i)F +A + M$
generate simultaneous jets
of order $s_1, \ldots, s_p$ at arbitrary distinct points
$x_1, \ldots , x_p$ of  $X$.

\noindent
Moreover, $ K_X + (n + 1)F +A + M$ is very ample. 
\end{lm}

\noindent
{\em Proof.} Fix $\omega$  a hermitian metric on $X$ and $g$ a
hermitian metric on $A$ with positive curvature
$\T{h}{A} \geq \frac{3}{2}\e \omega$ for some $\e>0$.

\noindent
 Since $F$ is ample and  the linear system $|F|$ 
is free of base-points,  for every
index $i$ there are $n$ sections $\{\s_{ij}\}_{j=1}^n$ of $F$ 
such that their common zero locus is zero-dimensional at $x_i$. 

\noindent
Define a s.h.m. $h$ on $(pn+\sum_i{s_i}) F $ by first defining 
 metrics $h_i$ on $(n+ s_i)F$:
$$
h^{-1}_i:= \left[ \sum_{j=1}^n{|\s_{ij}^2|} \right]^{n+s_i} 
$$
and then by multiplying them together
$$
h: = \prod_{i=1}^p h_i.
$$
\noindent
Since $M$ is nef, one can choose a 
hermitian metric $l$ on it such that
$\T{l}{M} \geq - \frac{1}{2} \e \omega$. 

\noindent
Define a metric $H$ on  $(pn + \sum_{i=1}^p s_i)F +A + M$ by setting
$$
H:=h\otimes g \otimes l.
$$
We have that $\T{h_i}{(n+s_i)F}\geq 0$, $\forall i$ so that 
$\T{h}{(np + \sum s_i)F} \geq 0$. It follows that $\Theta_H \geq \e 
\omega$. 

\noindent
Since $g$ and $l$ are continuous, $\id(H)=\id (h)$.

\noindent
By virtue of \ci{dem94}, Lemma 5.6.b, we have that
$\id (H)_{x_i}=\id (h)_{x_i} \subseteq {\frak m}_{x_i}^{s_i +1}$, 
$\forall i$
and that the scheme associated with $\id (H)$ is zero-dimensional
at all the points $x_i$.
We conclude by Proposition \ref{cucu}.

\noindent
The second part of the statement is  \ci{an-siu}, Lemma 11.1 
(the proof of which contains minor
inaccuracies but it is correct).
\blacksquare

\subsection{Effective results on vector bundles}
\label{vecbl}
We now see how to ``transplant" the metrics of Proposition
\ref{heart} to vector bundles
and how to use the results of \S\ref{vt} 
to prove effective results for  the  vector
bundles of the form $\frak P$  as in  the Introduction.

\medskip
Let us remark
that the lower bounds on $m$ given in the various statements of the 
theorem 
that follows are only indicative. Any improvement of these bounds
in the line bundle case that  can be obtained using 
strictly positive singular metrics would  give an analogous improvement
in the vector bundle case; see 
\ci{siu94b}, Proposition 5.1 for example.

\medskip
Let $n$, $p$, $\{s_1, \ldots , s_p \}$ and the various $m_i$ be as
in section \S\ref{demsiu}. Assume that $E$ is a rank $r$ vector bundle
on $X$
and let $N:=\min \{n, r \}$.
\begin{rmk}
\label{more}
{\rm 
See Ex.  \ref{listnef} for examples of $N$-nef vector bundles. 
}
\end{rmk}

\begin{tm}
\label{effres}
Let $X$ be   a projective manifold of dimension
$n$, $E$ be   $N$-nef, $A$  and $L$  be  ample 
line bundles on $X$.

\medskip
\noindent
{ \rm (\ref{effres}.1)} If $m\geq m_1(n,p)$, then
the global sections of $K_X\otimes E \otimes (mL)$ separate arbitrary
$p$ distinct points of $X$;

\medskip
\noindent
{ \rm (\ref{effres}.$1'$)} if $m\geq \frac{1}{2}(n^2 +n +2)$, then
$K_X\otimes E \otimes (mL)$ 
is generated by its global sections;

\medskip
\noindent
{ \rm (\ref{effres}.2)} if $m\geq m_2(n,p; s_1, \ldots, s_p)$, then
the global sections of $2K_X \otimes E \otimes (mL)$
generate simultaneous jets of order
$s_1,\ldots, s_p \in {\Bbb N}$ at arbitrary  $p$ distinct points of $X$;

\medskip
\noindent
{ \rm (\ref{effres}.$2'$)} if
$m\geq m_2(n,1;1)$, then 
the global sections of $2K_X \otimes E \otimes (mL)$
separate arbitrary pair of points of $X$ and
generate   jets of order
$1$  at an arbitrary  point of $X$.

\medskip
\noindent
{ \rm (\ref{effres}.3)} if $m\geq  m_3(n,p; s_1, \ldots , s_p)$, then
the global sections of 
$
(pn+\sum s_i  +1)K_X \otimes E \otimes (mL) \otimes A
$
generate simultaneous jets of order
$s_1,\ldots, s_p \in {\Bbb N}$ at arbitrary  $p$ distinct points of $X$;

\medskip
\noindent
{ \rm (\ref{effres}.4)} if
$m\geq  m_4(n)$, then 
the global sections of $(n+2)K_X \otimes E \otimes (mL) \otimes A$
separate arbitrary pair of points of $X$ and
generate   jets of order
$1$  at an arbitrary  point of $X$;

\medskip
\noindent
{ \rm (\ref{effres}.5)}
the global sections of
$E\otimes (mL)$  separate arbitrary pair of points of $X$ and
generate   jets of order
$1$  at an arbitrary  point of $X$ as soon as 
$$
m\geq C_n (L^n)^{3^{(n-2)}} ( n+2  + \frac{L^{n-1} \cdot 
K_X}{L^n})^{3^{(n-2)}(
\frac{n}{2} + \frac{3}{4})+ \frac{1}{4}},
$$
where $C_n=(2n)^{\frac{3^{(n-1)} - 1}{2}}(n^3-n^2 -n -1)^{
3^{(n-2)}(
\frac{n}{2} + \frac{3}{4})+ \frac{1}{4}}$.

\end{tm}

\begin{rmk}
\label{geoint}
{\rm 
Let us give a geometric interpretation to, say,  (\ref{effres}.$2'$). 
We employ the notation of \S\ref{demsiu}. Let $(X,E,L,m)$ be as in 
(\ref{effres}.$2'$). 
Let $E':= E \otimes (K_X + mL)$ and $L'=(r+1)(2K_X+ mL) + \det E$; note that
$h^0:=h^0(X,E')=\chi (X,E')$ and  that $L'$ is very ample.
Then there is a closed embedding 
$$
\phi := f \times g : X \longrightarrow G(r,h^0) \times \pn{n}
$$
such that   $E \simeq f^*(\QQ \otimes \q) \otimes g^*\odixl{\pn{n}}{-1}$,
$\deg \hat{f}(X)= (\det E')^n$ and $g$ is finite surjective with $\deg g 
=L'^n$. 

\noindent
Let $\{X_i,E_i,L_i\}_{i\in I}$ be a set 
 of triplets as above. If we can bound from
above $h^0_i$, $\deg \hat{f_i}(X_i)$ and ${L'_i}^n$, then we can find 
embeddings
$\phi_i : X_i \to G \times \pn{n}$ 
with 
$G=G(r,\max_I (h^0_i))$ such that the relevant invariants are bounded 
from above.
This applies, for example, to the set of flat
vector bundles of fixed rank over a (family of)
projective manifold(s), to the set of
all nef vector bundles of fixed rank over  curves of fixed genus,  
 to the set of projective surfaces with nef
 tangent bundles, etc. By virtue of Remark \ref{algnef},  a similar
remark holds, more generally, for nef vector bundles.
}
\end{rmk}

\medskip
\noindent
{\em Proof of Theorem \ref{effres}.} Note that 
$(5.2.2.1')$ and $(5.2.2.2')$ are special cases of $(5.2.2.1)$ and
$(5.2.2.2)$, respectively. We shall prove $(5.2.2.1)$ and $(5.2.2.2)$
in detail to illustrate the method. 
The remaining three assertions are left to the reader
and  can be proved using 
the same method with the aid of  Lemma \ref{freetojet} 
for the second and third to last,
and with the guideline of \ci{dem96}, 4.7 for the last one.

\medskip
\noindent
Proof of (\ref{effres}.1). We follow closely
\ci{an-siu}. The proof is by induction on $p$. Let $p=1$.
Let $x\in X$ be arbitrary. By (\ref{heart}.1) we have
a strictly 
positive s.h.m. $h$ on $mL$ such that 
$x$ is an isolated point of 
the scheme associated with
$\id (h)$.
By virtue of Theorem \ref{myvan},
 $H^1(X, K_X \otimes E \otimes (mL) \otimes \id (h))=0$ and 
the following
surjections imply the case $p=1$:
$$
H^0(X, K_X \otimes E \otimes (mL)) \surj H^0(X, K_X \otimes E \otimes 
(mL) 
\otimes
\odix{X}/\id (h))
$$
$$
 \surj H^0(X, K_X \otimes E \otimes (mL) \otimes
\odix{X}/ {\frak m}_x). 
$$
Let us assume that (\ref{effres}.1) is true for all integers
  $\rho \leq p-1$ and prove the case $\rho=p$. 
Let $h$ be as in (\ref{heart}.1) and  $\id (h)$ be  its  
multiplier ideal.  By  virtue of 
 Theorem \ref{myvan},
we have that  $H^1(K_X\otimes E \otimes (mL) \otimes \id (h))=0$.
Let $\cal J$ be the ideal sheaf on $X$ which agrees with
$\id (h)$ on $X\setminus J_0$ and which agrees with $\odix{X}$ on $J_0$.
Relabel the points so that $J_0=\{ 1, \ldots , l\}$.
By tensoring the exact sequence
$$
0 \to \id (h) \to {\cal J} \to {\cal J}/\id(h) \to 0
$$
with $K_X \otimes E \otimes mL$ 
we get the surjection: 
$$
H^0(X, K_X\otimes E \otimes (mL) \otimes {\cal J}) \surj 
\bigoplus_{i=1}^{l}   {\cal O}( K_X\otimes E \otimes (mL))_{x_i}
\otimes
\odix{X,x_{i}}/{\frak m}_{x_{i}}
$$
which implies that we can choose sections
${a_{1,j}} \in H^0(X, K_X \otimes E \otimes
(mL))$ vanishing at $x_2, \ldots , x_p$, but generating the stalk
$ (K_X \otimes E \otimes (mL))_{x_1}$. We now apply the induction
 hypothesis
to the set of $p-1$ points $\{x_2, \ldots , x_p\}$. \
By repeating the   above procedure, and keeping in mind
that at each stage we may have 
to relabel the points,
 we obtain  sections
$\{a_{i,j_{i}}\}\in H^0(X, K_X \otimes E \otimes
(mL))$, $\forall \, 1\leq i \leq p$ vanishing at $\{x_{i +1}, 
\ldots , x_p\}$
 but generating the stalk
$ (K_X \otimes E \otimes (mL))_{x_{i}}$.
Given any point $x_{i}$, with $1\leq i \leq r$, and any vector $w \in
 (K_X \otimes E \otimes (mL) )_{x_{i}}\otimes \odix{X, x_{i}}/
m_{x_{i}}$ 
it is now easy to find
 a linear combination of the sections $a_{i,j_i}$ which is $w$ at $x_{i}$
and zero at all the other $p-1$ points. This proves (\ref{effres}.1).

\medskip
\noindent
Proof of
(\ref{effres}.2).  
 We  fix  the integers $p,$
$s_1, \ldots, s_p$ and $p$ arbitrary distinct points
on $X$. We take a singular metric
$h$ on $K_X +mL$ with $m\geq m_1$ for which the associated multiplier
ideal $\id (h)$ has  the properties 
ensured by  (\ref{heart}.2). Theorem \ref{myvan} gives us the vanishing
of $H^1(K_X \otimes K_X \otimes E \otimes (mL) \otimes
\id(h))$ which, in turn, 
gives the
wanted surjection in view of the obvious surjections
$$
\odix{X,x_{i}}/\id (h)_{x_{i}} \to
 \odix{X,x_{i}}/{\frak m}_{x_{i}}^{s_{i} +1},
 \quad 
\forall \, \,  1 \leq i \leq p.
$$ 
\blacksquare

\begin{rmk}
{\rm 
Both statements in {\rm Proposition \ref{heart}}
 have counterparts  entailing
not powers  $mL$ of an ample line bundle $L$, but directly
an ample line bundle $\LL$ which has ``intersection theory"
large enough. See {\rm \ci{an-siu}}, {\rm Theorem 0.3}, 
{\rm \ci{dem96} Theorem
2.4.b}  and {\rm \ci{siu94b}}. 

\noindent
As a consequence one has statements similar to the ones
of {\rm Theorem \ref{effres}} with $mL$ substituted by an ample
line bundle $\LL$ with intersection theory large enough; we omit the details.
}
\end{rmk}

\begin{rmk}
\label{algnef}
{\rm 
Let $X$, $n$, $E$, $L$  be as in this section except 
that $E$ is only assumed
to be nef. Using algebraic techniques
we can see that  
 $K_X \otimes E \otimes \det \, E \otimes L^{m}$
is globally generated for $m\geq \frac{1}{2}(n^2 +n+2) $.
 Statements involving higher jets can be proved as well.
Details will appear in \ci{deeff}.
}
\end{rmk}

\begin{??}
{\rm
Let $X$, $E$ and $L$ be as above. Is the vector bundle
$K_X \otimes E  \otimes L^{\otimes m}$
generated by global sections for every $m\geq \frac{1}{2}(n^2 +n+2)$?
}
\end{??}

\bigskip
1991 {\em Mathematics Subject Classification}. 14C20, 14C30, 14F05, 
14F17, 14H60,
14M25, 14Q20, 32C17, 32C30, 32F05, 32J25, 32L10, 32L20, 53C55.

{\em Key words and phrases}. Singular hermitian metric, vector bundle,
transcendental methods,
curvature current, positivity, semi-positivity,
vanishing theorems, effective global generation, jets.

\bigskip
AUTHOR'S ADDRESS

\medskip
\noindent
Mark Andrea de Cataldo,
Max-Planck-Institut f\"ur Mathematik, Gottfried-Claren-Strasse
26, 53225 Bonn, Germany.

\noindent
e-mail: {\em markan@mpim-bonn.mpg.de}

\end{document}